\RequirePackage{fix-cm}
\documentclass[trsc,nonblindrev]{informs4}
\renewcommand{\theARTICLETOPLEFT}{\includegraphics{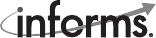}}
\renewcommand{\theARTICLETOPRIGHT}{}
\RRHSecondLine{}
\LRHSecondLine{}

\usepackage{tgtermes,newtxtext}
\DeclareFontFamily{OMS}{ntxtlf}{}
\DeclareFontShape{OMS}{ntxtlf}{m}{n}{<->ssub*cmsy/m/n}{}
\usepackage{newtxmath,bbm}
\input{ubbm.fd}
\DeclareFontShape{U}{bbm}{m}{n}{<-5.5>bbm5 <5.5-6.5>bbm6 <6.5-7.5>bbm7 <7.5-8.5>bbm8 <8.5-9.5>bbm9 <9.5-11>bbm10 <11-15>bbm12 <15->bbm17}{}
\usepackage{natbib,xcolor,xurl,hyperref,placeins,needspace}
\hypersetup{colorlinks=true,linkcolor=blue,citecolor=blue,urlcolor=blue,
  pdftitle={Learned spatial opportunity cost admission for ride-hailing with calibrated thresholds},
  pdfauthor={Baichuan Mo, Yuan Zhou, Qixin Zhang, Hai Yang}}
\bibpunct[, ]{(}{)}{,}{a}{}{,}
\def\bibfont{\small}
\def\bibsep{\smallskipamount}

\makeatletter
\newcommand{\settextleading}[2]{\fontsize{#1}{#2}\selectfont\@bls=\baselineskip}
\g@addto@macro\normalsizeXIInoCE{\settextleading{12}{19.32}}
\g@addto@macro\smallXnoCE{\settextleading{10}{16.1}}
\g@addto@macro\footnotesizeIXnoCE{\settextleading{9}{14.49}}
\makeatother
\SingleSpacedXII
\ECRepeatTheorems
\MANUSCRIPTNO{}
\EquationsNumberedThrough
\TheoremsNumberedThrough
\let\INFORMSproof\proof
\let\endINFORMSproof\endproof
\renewenvironment{proof}[1][Proof]{\INFORMSproof{#1.}}{\endINFORMSproof}

\begin{document}
\RUNAUTHOR{Mo, Zhou, Zhang, and Yang}
\RUNTITLE{LSOA for ride-hailing with calibrated thresholds}
\TITLE{Learned spatial opportunity cost admission for ride-hailing with calibrated thresholds}
\ARTICLEAUTHORS{%
  \AUTHOR{\href{mailto:bmo@tsinghua.edu.cn}{\textcolor{black}{Baichuan Mo}}\textsuperscript{1}, \href{mailto:yuan-zhou@tsinghua.edu.cn}{\textcolor{black}{Yuan Zhou}}\textsuperscript{2,*}, \href{mailto:qixin.zhang2026@gmail.com}{\textcolor{black}{Qixin Zhang}}\textsuperscript{3}, \href{mailto:cehyang@ust.hk}{\textcolor{black}{Hai Yang}}\textsuperscript{4}}%
  \AFF{\textsuperscript{1}Department of Civil Engineering, Tsinghua University, Beijing, China.\\
    \textsuperscript{2}Yau Mathematical Sciences Center, Tsinghua University, Beijing, China.\\
    \textsuperscript{3}Wuhan University, Wuhan, China.\\
    \textsuperscript{4}Department of Civil and Environmental Engineering, The Hong Kong University of Science and Technology, Hong Kong, China.\\
    \textsuperscript{*}Corresponding author. \href{mailto:yuan-zhou@tsinghua.edu.cn}{\textcolor{black}{yuan-zhou@tsinghua.edu.cn}}}%
}
\ABSTRACT{Ride-hailing operators must accept or reject requests before later demand is observed. Serving a ride earns revenue but occupies a vehicle and changes where it can serve next. We study immediate admission for a centrally managed fixed fleet and develop learned spatial opportunity cost admission. The policy estimates demand from historical and observed requests and solves a rolling mean-demand fleet plan. Its dual prices charge for origin capacity and credit availability at the destination after completion. With positive regularization, admission is equivalent to applying a threshold to a unique planned service fraction. We hold regularization fixed and calibrate the threshold using realized validation revenue. We characterize when a common threshold attains optimal revenue and bound the loss when it cannot implement the required dispatch decisions. We obtain exact optimality conditions on a single route. When several origins share return demand at a hub, the revenue loss under the best common threshold can grow in proportion to fleet size. Independent validation provides a revenue guarantee relative to the best candidate threshold. Under Poisson arrivals and a margin condition, the revenue difference between the same policy using estimated and true demand forecasts decreases at a tight inverse-square-root rate as historical data accumulates. Experiments cover 60 synthetic instances in six demand families and ten New York workdays across 25 taxi zones and three fleet sizes. With two validation days, calibration yields mean paired revenue gains of 4.93\% in the synthetic study and 7.08\% in New York relative to the same policy with a zero adjusted-reward cutoff. In New York, mean paired revenue is 1.24--2.54\% below the best-performing benchmark, while average decision time is about one-eleventh as large. Calibration gains increase as fleet size falls in New York. In the synthetic study, the largest gains over accepting every feasible request occur under commuting demand. }
\maketitle

% Revision author: Baichuan Mo
\section{Introduction}\label{sec:introduction}

Ride-hailing operators must decide which requests to serve with a limited fleet. Accepting a ride earns revenue, but the vehicle is unavailable for other requests until the trip ends. A long trip may prevent it from serving several later requests, whereas a trip to a region with high demand may allow it to serve another passenger soon after drop-off. The operator therefore needs to consider both the fare and the rides the vehicle could serve afterward. Here we study sequential admission for an operator with a centrally managed fleet: it accepts or rejects each request when it arrives, before observing later requests. Trip rewards and fleet participation are given. For example, a chauffeur company working with Blacklane (a chauffeur service platform) sees the trip details and the payment it would receive before deciding whether to accept the ride \citep{blacklane2026partners}. On Booking.com, a partner fleet quotes a price for the trip. If the passenger books and Booking.com assigns the trip to that fleet, the fleet decides whether to accept the booking at the quoted price \citep{bookingTaxiSupplier}. The key research question is how can the operator estimate the opportunity cost of committing a vehicle in real time, accounting for its future location and release time while using imperfect demand forecasts?

\begin{figure}[htbp]
\centering
\includegraphics[width=0.93\textwidth]{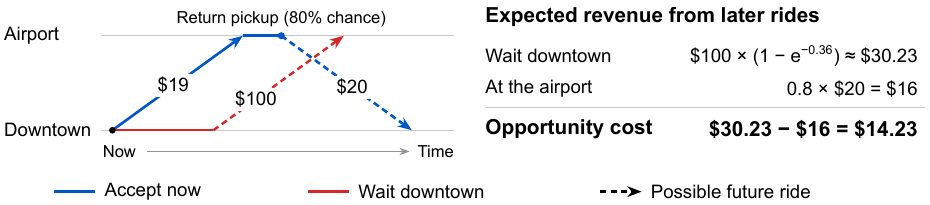}
\caption{An example illustrating the opportunity cost of dispatching a vehicle to the airport.}
\label{fig:airport-decision}
\end{figure}

Figure~\ref{fig:airport-decision} illustrates this opportunity cost. One vehicle is idle downtown when a passenger offers \$19 for a ride to the airport. A \$20 return request arrives at the airport with probability 0.8, and only the current ride gets there in time to serve it. If the vehicle waits downtown, it can serve at most one later ride to the airport. Each such ride pays \$100, and the number of requests follows a Poisson distribution with mean 0.36. Waiting therefore earns $100(1-e^{-0.36})\approx\$30.23$ in expectation. Sending the vehicle now gives up that opportunity but earns an expected $0.8\times\$20=\$16$ from the return ride. The opportunity cost of dispatch is thus $\$30.23-\$16=\$14.23$. The current reward of \$19 covers this cost, so the operator should accept.

A general way to express this opportunity cost is to compare future revenue after accepting and rejecting the current request. For a physically feasible request $i$ of type $\xi_i$ and reward $r_{\xi_i}$, let $\boldsymbol s_i^{\mathrm{accept}}$ and $\boldsymbol s_i^{\mathrm{reject}}$ denote the post-decision states after acceptance and rejection, respectively. Each state records vehicle locations and release times, the current time, and the observed demand history. Let $V_i^*(\boldsymbol s)$ be the optimal expected future revenue from state $\boldsymbol s$ under the true demand model, excluding the current reward and using only information available as requests arrive. The opportunity cost of dispatching a vehicle for this request is $c_{\mathrm{opp},i}^*:=V_i^*(\boldsymbol s_i^{\mathrm{reject}})-V_i^*(\boldsymbol s_i^{\mathrm{accept}})$. The operator should accept when $r_{\xi_i}\ge c_{\mathrm{opp},i}^*$. Computing this cost requires accounting for how dispatch changes the fleet's subsequent service opportunities.

Existing methods use approximate value functions, rollout simulation, or deterministic plans based on expected demand \citep{li2023bidprices,bertsekas1999rollout,bumpensanti2020resolving}. Exact continuation values are difficult to compute because future admission decisions generate a large fleet-state space. Simulation-based approximations avoid this full recursion but require repeated evaluation of candidate actions. Replacing future demand by its mean reduces computation, but errors in the two continuation values can overstate or understate their difference and change the admission decision. In the airport example, the mean-demand plan values waiting downtown at $0.36\times\$100=\$36$. This overstates what the vehicle can earn. After subtracting the expected \$16 from the return ride, the plan estimates the opportunity cost as $\hat c_{\mathrm{opp},i}=\$36-\$16=\$20$. The \$19 reward of the current ride is below this estimate, so the plan rejects the ride even though accepting it is optimal. Hence, how can the operator use mean demand to estimate opportunity costs while limiting the effect of these approximation errors on which requests it accepts?

One idea is to use the mean-demand estimate $\hat c_{\mathrm{opp},i}$ as a starting point and adjust the minimum reward required to accept a ride. In the airport example, lowering that cutoff from \$20 to \$19 would lead the operator to accept the current ride, which is the optimal decision. This example motivates calibrating the threshold for accepting ride-hailing requests. Prior work has studied threshold rules for allocating a fixed resource budget \citep{bumpensanti2020resolving,vera2020prophet}. However, it remains unclear when a single calibrated threshold can recover optimal admission decisions across a fleet and, when it cannot, how much revenue the fleet loses even with the best threshold.

We propose learned spatial opportunity cost admission (LSOA), a policy that estimates opportunity costs from a rolling fleet plan based on mean demand and uses validation data to calibrate the threshold for accepting a request. We aim to answer three questions. First, suppose the operator knows the true mean demand. When can it earn the same expected revenue as the optimal policy by using one threshold for all requests? If no single threshold achieves this, how much more can the optimal policy earn than the best threshold policy? Second, if the operator selects a threshold using a limited number of validation days, how close is its expected revenue to that of the best candidate threshold? Third, when mean demand must be estimated, which admission decisions can change, and how quickly does the revenue difference from using the true mean demand shrink as more historical data become available?

Our analysis first shows when a common threshold attains optimal revenue. When only some vehicles are needed for return trips at a hub, a common threshold can send all vehicles or none if the plan assigns the same fraction to every outbound request. The resulting revenue loss can grow in proportion to fleet size. With $D^{\mathrm{val}}$ validation days, the revenue gap to the best candidate threshold is $O((D^{\mathrm{val}})^{-1/2})$. At a fixed threshold, small forecast errors can change decisions when the planned service fraction is close to the cutoff. With $D$ historical days, the resulting revenue gap relative to using true mean demand has a tight $D^{-1/2}$ bound.

We evaluate LSOA on 60 synthetic instances spanning six demand families and a controlled replay of ten New York workdays in 25 Manhattan taxi zones at three fleet sizes. In the synthetic study, LSOA earns higher mean revenue than the Greedy benchmark in every demand family and has the highest overall average revenue among the eight online methods. Its mean paired gain over Greedy is 7.90\%, with larger gains under spatially concentrated and commuting demand. In New York, these gains range from 6.53\% to 20.56\%, while threshold calibration is most valuable at the smallest fleet. Across the three fleet sizes, LSOA earns 1.24\%--2.54\% less revenue than the best-performing benchmark in mean paired comparisons, while requiring about one-eleventh of its average decision time. Together, these results show a consistent improvement over Greedy admission and a practical balance between revenue and computation.

The paper makes three contributions.
\begin{enumerate}
\item We develop a fleet admission policy that calibrates a common threshold on the fraction of a vehicle assigned to each incoming request by a rolling mean-demand plan. We derive the optimal threshold interval on a single route (Proposition~\ref{prop:calibration-value}) and quantify the revenue loss of the best common threshold when several origins share return demand at a hub (Proposition~\ref{prop:hub-threshold}).

\item We bound the revenue effects of selecting a threshold and estimating demand. With $D^{\mathrm{val}}$ independent validation days, the expected revenue of the selected threshold is within $O((D^{\mathrm{val}})^{-1/2})$ of the best candidate in a fixed set with high probability. At a fixed threshold, we bound the change in the planned service fraction by the total error in the forecast request counts. This identifies which admissions can change and yields a tight $D^{-1/2}$ bound on the revenue gap relative to using true mean demand with $D$ historical days under Poisson arrivals and a margin condition.

\item We quantify the value of threshold calibration and identify when accounting for future requests is most useful. Using synthetic demand and a replay of New York taxi requests, we examine how these gains vary with fleet size and demand patterns and compare policies in terms of revenue and computation time.
\end{enumerate}

We summarize how our setting and admission policy differ from those in the closest studies in Section~\ref{sec:key-references}.

\section{Literature review}\label{sec:literature}

The paper connects three streams: admission with spatially reusable capacity, planning approximations for online control, and learning with revenue guarantees.

\subsection{Admission with spatially reusable capacity}

Serving a request temporarily occupies a resource. Reusable-resource models capture this feature in ride-sharing matching \citep{dickerson2018allocation}, assortment optimization \citep{gong2022assortment}, and allocation with stochastic usage times \citep{goyal2025reusable}. For vehicles, reuse also changes location. \citet{li2023bidprices} study this coupling through Lagrangian bid-price controls for a car rental network with booking lead times and rental durations. \citet{kanoria2022mirror} develop supply-dependent admission prices for circulating capacity without estimating demand rates.

Related ride-hailing studies examine spatial pricing \citep{bimpikis2019spatial,chen2024spatial}, empty-car routing \citep{braverman2019empty}, and admission and repositioning with strategic drivers \citep{afeche2023networks}. Matching policies use continuous linear programs \citep{ozkan2020matching}, account for pickup travel \citep{feng2021way}, or combine learned vehicle values with assignment \citep{xu2018dispatch}. These studies establish the importance of future vehicle availability in current decisions.

Our decision is immediate acceptance or rejection with given fares and a centrally managed fixed fleet, rather than pricing or driver incentives. LSOA obtains prices from a rolling plan that includes the current request, forecast services, and scheduled releases of busy vehicles. The origin charge and destination credit account for both where and when capacity can be reused. The admission cutoff is then calibrated by realized policy revenue, rather than fixed by the planning prices alone.

\subsection{Planning approximations and admission thresholds}

Bid-price control converts a resource valuation into an acceptance test. \citet{talluri1998bidprice} show that such controls need not be optimal at finite scale, but establish asymptotic optimality when capacity and sales volumes grow. \citet{talluri1999randomized} compute prices using linear programs on sampled demand. Re-solving updates deterministic plans as arrivals occur, with revenue-loss guarantees in network revenue management \citep{jasin2012resolving,bumpensanti2020resolving}. \citet{brown2025fluid} analyze feasible policies based on fluid relaxations, including reoptimized value functions and adjustments that respect shared resource constraints.

Continuation values provide another route to opportunity costs. \citet{topaloglu2007sensitivity} evaluate marginal fleet capacity through approximate dynamic programming, and rollout evaluates candidate actions using a continuation policy \citep{bertsekas1999rollout}. Sampling demand and optimizing each complete path is a different approximation: those optimizations know future requests before deciding. Information-relaxation theory formalizes this distinction between advance information and feasible online decisions \citep{brown2010information}.

LSOA combines tractable mean-demand planning with revenue-based calibration of the admission rule. Rather than refining continuation-value approximations or applying planning prices with a fixed cutoff, we calibrate how approximate prices guide accept-or-reject decisions using realized revenue from independent validation request streams. This targets the decision impact of planning-model distortion, rather than demand forecasting alone.

\subsection{Learning, calibration, and revenue guarantees}

Online allocation methods combine resource control with learning. \citet{balseiro2023mirror} analyze dual mirror descent across several input models, and \citet{zhang2025reusable} combine dual learning with tests for nonstationarity in reusable-resource allocation. Their guarantees address policy performance under uncertain inputs; the latter approaches optimal reward as resource inventories and the horizon increase.

Our focus is revenue-based selection of the admission threshold for a given regularization parameter. We distinguish the gain from selecting a better admission rule from the revenue effect of learning demand, providing separate calibration and tight $D^{-1/2}$ learning guarantees.

\subsection{Differences from key references}\label{sec:key-references}

\citet{li2023bidprices} study whether to accept advance and on-the-spot car rental bookings. Their bid prices account for both the pickup station and the return station, where the car becomes available the next day. In ride-hailing, a vehicle leaves as soon as the operator accepts a request and can serve another passenger as soon as the trip ends. Even two rides on the same route can therefore have different opportunity costs if only the earlier ride reaches the destination in time for a later pickup. We estimate this cost using a rolling plan based on mean demand (rather than a Lagrangian relaxation of a stochastic dynamic program). The plan tracks when each busy vehicle becomes available, but mean demand can overstate the reward from waiting for a pickup that may never arrive. This gives a reason to calibrate the cutoff for dispatch even when demand means are known.

In the ride-hailing model of \citet{kanoria2022mirror}, the platform can reject a request or choose among eligible vehicle locations from which to serve it. Their main analysis assumes instantaneous trips. With positive travel times, their mirror backpressure rule uses current idle supply and a charge for time spent on the trip. In our setting, a ride requires an idle vehicle at its origin, and accepting it commits that vehicle until drop-off. Our rolling plan uses mean demand and scheduled vehicle returns to estimate the value of a vehicle at the origin now and at the destination after drop-off. The difference gives an estimated opportunity cost for each ride. To account for errors in that estimate, we use validation revenue to choose one adjustment that raises or lowers the reward required for acceptance across all rides. 

\citet{akan2009bidprice} study a network with fixed capacities in which the operator chooses how much of the demand rate for each product to accept. Their quadratic perturbation yields a unique booking rate, which the operator can implement directly in the rate-based model. A ride request, however, requires a whole vehicle. In our plan, the quadratic term gives a unique planned fraction for that request. Its marginal value links a threshold on the fraction to a cutoff on the trip reward net of the estimated opportunity cost. This link also lets us study forecast error. A change in the planned fraction changes the dispatch decision only if it moves the fraction across the threshold.

In online packing with fixed resource budgets, \citet{vera2020prophet} use a cutoff of one-half on the planned share of demand served in the arriving class. \citet{bumpensanti2020resolving} round class acceptance probabilities near zero or one and use the planned probabilities between these cutoffs. In our fleet, accepting a ride also makes a vehicle available for later pickups at its destination. When several origins supply the same hub, the value of sending one vehicle can depend on how many others have already been sent there. We therefore ask when the fleet can make optimal dispatch decisions using the same threshold on planned service fractions at every origin.

In \citet{balseiro2021relocating}, the provider quotes a price for each request, and a customer decides whether to pay it. In our setting, the reward is given and the fleet operator decides whether to accept the ride. 

\citet{topaloglu2008stochastic} studies an airline that accepts a feasible itinerary when its revenue exceeds the sum of the bid prices of the flight legs it uses. Simulated bookings are used to optimize those prices. In our setting, the cost of dispatch depends on which other vehicles will be available at the origin and destination when they are needed. We therefore recompute the origin charge and destination credit from the current fleet plan at each feasible request. With regularization fixed, only the common cutoff on the adjusted reward is calibrated. We compare candidate cutoffs by replaying each one over an entire validation day, because changing one dispatch can change which later rides the fleet can serve.

\section{Problem description}\label{sec:problem}

We study how a ride-hailing operator with a fixed fleet decides which requests to accept. Each request is considered when it arrives. Historical orders help estimate future demand, and vehicle records show where and when capacity will become available.

\subsection{Requests and time intervals}

The service area has $S$ regions indexed by $s\in\mathcal S$, where $\mathcal S:=\{1,\ldots,S\}$. The fleet has $C$ vehicles indexed by $c\in\mathcal C$, where $\mathcal C:=\{1,\ldots,C\}$. We divide the operating horizon $[0,H)$ into $T$ time intervals of equal length $\Delta$, so that $H=T\Delta$. The intervals are indexed by $\mathcal T:=\{0,\ldots,T-1\}$, with interval $t\in\mathcal T$ covering $[t\Delta,(t+1)\Delta)$. This time grid is used to aggregate demand and plan vehicle services.

Requests are indexed separately in arrival order by $i\in\mathcal I$, where $\mathcal I:=\{1,\ldots,I\}$. The total number $I$ is random. Request $i$ arrives at clock time $t_i\in[0,H)$, with $t_1\le\cdots\le t_I$. A fixed priority order breaks timestamp ties. Thus $i$ identifies a request, $t_i$ is its arrival time, and $t$ identifies a time interval. Request $i$ belongs to interval $t$ when $t\Delta\le t_i<(t+1)\Delta$. Several requests can arrive within the same interval or share an arrival timestamp.

There are $Z$ known request types indexed by $z\in\mathcal Z$, where $|\mathcal Z|=Z$. For example, requests can be classified into types by their origin--destination pairs. Then we have $\mathcal Z=\mathcal S\times\mathcal S$ and $Z=S^2$. Type $z$ has origin $o_z$, destination $d_z$, total trip reward $r_z\ge0$, and duration $\ell_z\Delta$, where $\ell_z$ is a positive integer. Let $\xi_i\in\mathcal Z$ be the random variable giving the type of request $i$. Its arrival time $t_i$ and type $\xi_i$ are observed before its decision. An interval without a request contains no request index. For example, if $\Delta=2$ minutes, arrivals at 08:00:10, 08:00:45, and 08:01:30 are three requests in one two-minute interval. They receive three decisions sequentially.

\subsection{Vehicle state and admission}

Immediately before request $i$, $s_{c,i}$ and $u_{c,i}$ record the next available region and clock time of vehicle $c$. The subscript $i$ indicates which request is being considered. Vehicle $c$ is idle if $u_{c,i}\le t_i$. Otherwise, the vehicle remains busy until $u_{c,i}$. Initially busy vehicles finish their existing trips as scheduled, and the rewards of those existing trips are excluded from the planning objective.

The binary decision $a_i$ equals one if request $i$ is accepted and zero if it is rejected. Acceptance requires an idle vehicle at $o_{\xi_i}$ and $t_i+\ell_{\xi_i}\Delta\le H$. It earns $r_{\xi_i}$ and makes the assigned vehicle next available at $d_{\xi_i}$ at clock time $t_i+\ell_{\xi_i}\Delta$. Rejection leaves the fleet unchanged and the request leaves. An infeasible request is rejected.

Vehicles can pick up immediately within their region, serve one request at a time, and move only while serving a ride. The fleet size $C$ stays fixed. The model excludes empty relocation, passenger queues, and random travel times. A sequential decision is made at each arrival. 

\subsection{Demand information and objective}

Demand observations consist of historical request streams and the current stream up to the decision being made. The platform records all arrivals, including rejected and infeasible requests. It also records observation durations, so periods without orders contribute time exposure.

Let $\mathbb P$ denote the joint probability distribution of the historical data and the operating-day stream $(I,(t_i,\xi_i)_{i=1}^I)$. This distribution may include variation in arrival rates and dependence between requests. We assume $I$ is finite almost surely and $\mathbb E_{\mathbb P}[I]<\infty$. Arrivals are exogenous to acceptance decisions. A feasible policy $\pi$ selects decision $a_i^\pi$ using the observed data, current request, and fleet records. Its expected revenue is
\begin{equation}\label{eq:policy-objective}
V_\pi(\mathbb P):= \mathbb E_{\mathbb P}\!\left[\sum_{i=1}^{I}r_{\xi_i}a_i^\pi\right].
\end{equation}
The expectation is over the historical sample, the number of requests $I$, their arrival times $t_i$, and types $\xi_i$. For randomized policies, the expectation also includes their internal randomness. Although historical data are observed before operations, the expectation averages over their possible realizations to evaluate the learning procedure across historical samples. We assume the fleet size is fixed, all accepted rides finish by $H$, and vehicles have zero terminal value.

\section{Demand estimation and fleet planning}
\label{sec:methodology}\label{sec:planning}

When a physically feasible request arrives, the algorithm first forecasts future request counts using historical data and arrivals observed so far. It then plans vehicle use over a moving time window. The plan accounts for the locations of idle vehicles and where and when busy vehicles will finish their trips. From this plan, the algorithm obtains spatial dual prices to evaluate the current request. The notation $|i$ means that a quantity is calculated for request $i$. Several requests may arrive within the same time interval $t$.

\subsection{Arrival-rate estimation}\label{sec:forecast-accuracy}

Clock time is partitioned into $J$ fixed demand contexts, such as hours of the day. These contexts are indexed by $j\in\mathcal J$, where $\mathcal J:=\{1,\ldots,J\}$. Different contexts may have different request arrival rates and origin--destination patterns. We therefore estimate the arrival rate of each request type separately for each context. The known function $\kappa:[0,H)\to\mathcal J$ maps each time $u$ to its demand context $\kappa(u)\in \mathcal J$. For each type and context, let $N^{\mathrm{hist}}_{z,j}$ be the historical request count and $E^{\mathrm{hist}}_j>0$ the historical observation duration. Before the admission decision for request $i$, the arrival-rate estimate for each type and context is updated using the available observations:
\begin{equation}\label{eq:demand-estimate}
\widehat\alpha_{z,j|i} :=\frac{N^{\mathrm{hist}}_{z,j} +\sum_{k<i}\mathbbm{1}\{\xi_k=z,\ \kappa(t_k)=j\}} {E^{\mathrm{hist}}_j+\int_0^{t_i}\mathbbm{1}\{\kappa(u)=j\}\,du}, \qquad \forall z\in\mathcal Z,\ j\in\mathcal J.
\end{equation}
Here $\mathbbm{1}\{\cdot\}$ denotes the indicator function. The resulting estimate $\widehat\alpha_{z,j|i}$ measures the arrival rate in requests per unit time. The numerator includes every previous request, regardless of its admission decision. The denominator includes all observation time, including periods with no arrivals. Request $i$ enters the count after its decision. Earlier requests with the same timestamp are included according to the fixed priority order.

For example, 120 historical requests over 60 minutes and 20 earlier requests over 10 minutes give a rate of 2 requests per minute. The expected count over the next 2 minutes is therefore 4 requests. Equation~\eqref{eq:demand-estimate} estimates a constant mean rate within each context. Rate changes and dependence can create forecast errors.

\subsection{The planning time grid}

The planning window covers up to $L$ intervals. When request $i$ arrives, the remaining operating time limits this window to $T_i:=\min\{L,\lceil(H-t_i)/\Delta\rceil\}$ intervals indexed by $\tau\in\mathcal T_i$, where $\mathcal T_i:=\{0,\ldots,T_i-1\}$. Planning interval $\tau$ covers $[t_i+\tau\Delta,\min\{H,t_i+(\tau+1)\Delta\})$. The grid moves with the current arrival, while each full interval still has length $\Delta$. For example, $\Delta=2$ minutes and $L=15$ give a 30-minute planning window.

The forecast number of additional type-$z$ requests in interval $\tau$ is
\begin{equation}\label{eq:planning-counts}
\widehat q_{z,\tau|i} :=\int_{t_i+\tau\Delta}^{\min\{H,t_i+(\tau+1)\Delta\}} \widehat\alpha_{z,\kappa(u)|i}\,du, \qquad \forall \tau\in\mathcal T_i,\ z\in\mathcal Z.
\end{equation}
These counts $\widehat q_{z,\tau|i}$ concern orders that have yet to be observed. In particular, $\widehat q_{z,0|i}$ describes additional demand in the current planning interval. The known request $i$ is represented separately by one unit.

The plan places forecast rides at interval starts. A busy vehicle enters the grid at the first planning time at or after its physical completion. A forecast ride is eligible if its planned completion is by $H$. If the ride finishes after the planning window, its full reward is still included, but the vehicle provides no further capacity within that window. These are planning formulation approximations. Actual accepted rides start at their request arrival times and release vehicles at their physical completion times.

\subsection{Vehicle service planning}

Let $\mathcal Z_{\tau|i}:=\{z\in\mathcal Z: t_i+(\tau+\ell_z)\Delta\le H,\ \widehat q_{z,\tau|i}>0\}$ contain the forecast types with positive expected counts whose rides, if started at planning interval $\tau$, would finish by the end of operations $H$. The decision variable $x_{z,\tau|i}$ is the planned number of additional type-$z$ rides accepted starting at that interval. The decision variable $y_i\in[0,1]$ is the planned service fraction of the current observed request $i$. Both variables are continuous in the planning problem, whereas the implemented admission $a_i$ is binary. The current-request fraction $y_i$ is kept distinct from the forecast amounts $x_{z,\tau|i}$ so that the model supports a decision at each arrival without representing every possible future request sequence. The forecast amounts $x_{z,\tau|i}$ guide the current decision but are not dispatch commitments. Only the current admission $a_i$ is implemented, subject to actual vehicle availability at $t_i$, using the rule in Section~\ref{sec:algorithm}. The plan is recomputed at the next feasible arrival using updated demand and vehicle records, so decisions need not wait for a planning boundary.

For $s\in\mathcal S$ and $\tau'\in\mathcal T_i$, define $b_{s,\tau'|i}:=\sum_{c\in\mathcal C} \mathbbm{1}\{s_{c,i}=s,\ u_{c,i}\le t_i+\tau'\Delta\}$. The supply $b_{s,\tau'|i}$ counts vehicles from the current records available at $s$ by that planning time. The regularization parameter $\theta\ge0$ is given before operations begin. For each feasible request $i$, LSOA determines the planned service amounts by solving the following problem:
\begin{subequations}\label{eq:quadratic-plan}
\begin{align}
\max_{y_i,\{x_{z,\tau|i}\}}\quad &F_i:=r_{\xi_i}y_i+\theta y_i(1-y_i) +\sum_{\tau\in\mathcal T_i}\sum_{z\in\mathcal Z_{\tau|i}} \left[r_zx_{z,\tau|i} +\theta x_{z,\tau|i}\left(1-\frac{x_{z,\tau|i}}{\widehat q_{z,\tau|i}}\right)\right],
\label{eq:quadratic-objective}\\
\text{s.t.}\quad &\mathbbm{1}\{o_{\xi_i}=s\}y_i +\sum_{\tau=0}^{\tau'}\sum_{z\in\mathcal Z_{\tau|i}:o_z=s}x_{z,\tau|i}\notag\\
&\quad\le b_{s,\tau'|i} +\mathbbm{1}\{d_{\xi_i}=s,\ \ell_{\xi_i}\le\tau'\}y_i +\sum_{\tau=0}^{\tau'}\sum_{\substack{z\in\mathcal Z_{\tau|i}\\d_z=s,\ \tau+\ell_z\le\tau'}} x_{z,\tau|i}, \quad \forall s\in\mathcal S,\ \tau'\in\mathcal T_i,
\label{eq:quadratic-vehicle}\\
&0\le y_i\le1,\qquad 0\le x_{z,\tau|i}\le\widehat q_{z,\tau|i}, \quad \forall \tau\in\mathcal T_i,\ z\in\mathcal Z_{\tau|i}.
\label{eq:quadratic-bounds}
\end{align}
\end{subequations}
When $\theta=0$, \eqref{eq:quadratic-plan} is a linear program (LP) that maximizes planned trip revenue. This is the linear special case of LSOA. For $\theta>0$, the same formulation adds quadratic regularization to the planned service fractions. Both cases use the same demand forecasts and vehicle constraints.

The left side of \eqref{eq:quadratic-vehicle} counts departures through interval $\tau'$. The right side of \eqref{eq:quadratic-vehicle} counts inherited capacity $b_{s,\tau'|i}$ and completed planned rides separately. The current ride consumes its planned fraction $y_i$ at its origin $o_{\xi_i}$ and returns that fraction at its destination $d_{\xi_i}$ after $\ell_{\xi_i}$ intervals. Forecast rides follow the same accounting. This keeps the planned idle inventory nonnegative in every region.

The objective in \eqref{eq:quadratic-objective} combines planned trip revenue with quadratic regularization. For a planned service amount $x$ with estimated bound $\widehat q$, this term is $\theta x(1-x/\widehat q)$. It is zero at $x=0$ and $x=\widehat q$ and positive between these endpoints when $\theta>0$. The ratio $x/\widehat q$ is the planned fraction of demand served. The current request has bound one and regularization term $\theta y_i(1-y_i)$.

For $\theta>0$, regularization makes the objective $F_i$ strictly concave in the service variables $y_i$ and $x_{z,\tau|i}$. It selects a unique service plan and gives curvature for the sensitivity analysis below. Its marginal contribution falls from $+\theta$ to $-\theta$ as the planned fraction increases from zero to one. The parameter has the same monetary units as a trip reward. Changing $\theta$ can change both the service allocation and the vehicle prices. We fix $\theta$ when calibrating the admission threshold.

The plan must respect both vehicle location and trip duration: capacity used for a ride can be reused only at its destination after completion. Proposition~\ref{prop:ts-duality} shows that the constraints enforce this requirement on the planning grid, even for fractional service amounts, and that $\theta>0$ gives a unique optimal service plan.

\begin{proposition}[Feasibility and uniqueness of the plan]
\label{prop:ts-duality}
For finite nonnegative forecasts and $\theta>0$, the feasible solutions of \eqref{eq:quadratic-plan} are fractional vehicle flows on the planning grid, and their optimal service amounts $y_i^*$ and $x_{z,\tau|i}^*$ exist and are unique.
\end{proposition}

\begin{proof}
See Section~\ref{ec:planning} of the e-companion.\quad$\square$
\end{proof}

Proposition~\ref{prop:ts-duality} has two practical implications. Feasibility means that the plan respects where and when vehicle capacity is available. Capacity assigned to a ride can be used again only at its destination after completion on the planning grid. Uniqueness means that a fixed fleet, forecast, planning grid, and $\theta>0$ give one optimal set of planned service amounts. In particular, there is a single current planned fraction $y_i^*$ to use in the admission rule.

\section{Dual prices and online admission}\label{sec:algorithm}

Solving \eqref{eq:quadratic-plan} yields dual multipliers for the vehicle constraints at each region $s$ and planning interval $\tau$. In this section, we use these multipliers to estimate the opportunity cost of serving the current request, and then update the physical fleet before considering the next arrival. This step connects the aggregate plan to sequential execution.

\subsection{Dual prices}\label{sec:vehicle-prices}

Problem~\eqref{eq:quadratic-plan} jointly accounts for forecast demand $\widehat q_{z,\tau|i}$ and scheduled vehicle availability $b_{s,\tau|i}$. It can be directly solved using an optimization solver, such as Gurobi. Along with the optimal plan, the solver returns the dual multiplier $\lambda_{s,\tau'|i}\ge0$ associated with the cumulative vehicle constraint \eqref{eq:quadratic-vehicle} for region $s$ and interval $\tau'$. This multiplier represents the marginal value of relaxing that constraint by one unit of vehicle capacity in the regularized objective, measured in monetary units per vehicle. If the constraint has slack at the optimum, its multiplier $\lambda_{s,\tau'|i}$ is zero.

A vehicle available at planning interval $\tau$ relaxes the cumulative constraints in its region from that interval onward. We define the dual price of that vehicle availability as the sum of their multipliers:
\begin{equation}\label{eq:vehicle-price}
\rho_{s,\tau|i}:=\sum_{\tau'=\tau}^{T_i-1}\lambda_{s,\tau'|i}, \qquad \forall s\in\mathcal S,\ \tau\in\mathcal T_i.
\end{equation}
For $\tau\ge T_i$, we set $\rho_{s,\tau|i}=0$, since capacity returned beyond the planning window has no terminal value in the plan.

For the observed request, the adjusted reward is $g_i:=r_{\xi_i}-\rho_{o_{\xi_i},0|i} +\rho_{d_{\xi_i},\ell_{\xi_i}|i}$. The origin term $\rho_{o_{\xi_i},0|i}$ accounts for using an available vehicle now. The destination term $\rho_{d_{\xi_i},\ell_{\xi_i}|i}$ accounts for restoring its availability at the destination after the trip. Thus, \(g_i\) equals the trip revenue \(r_{\xi_i}\) minus the estimated opportunity cost of serving the request, \(\rho_{o_{\xi_i},0|i}-\rho_{d_{\xi_i},\ell_{\xi_i}|i}\).

\subsection{Calibrated threshold admission}

For a given $\theta>0$, the admission threshold $\gamma\in(0,1)$ is selected before the operating day. Using the adjusted reward $g_i$ defined above, LSOA makes its decision for a feasible request $i$ according to
\begin{equation}\label{eq:price-rule}
a_i:=\mathbbm{1}\{g_i\ge\theta(2\gamma-1)\}.
\end{equation}
The platform holds $\theta$ fixed and replays each candidate threshold on validation request streams collected before the operating day. It selects the $\gamma$ with the highest mean total revenue and holds both parameters fixed during operations. The influence of $\theta$ and joint parameter calibration are examined in Section~\ref{ec:theta-sensitivity} of the e-companion. When $\theta=0$, the rule reduces to $a_i=\mathbbm{1}\{g_i\ge0\}$, independently of $\gamma$. Positive regularization provides the adjustable cutoff used for calibration. Sections~\ref{sec:calibration-value}--\ref{sec:calibration-selection} explain when threshold adjustment improves mean-demand decisions and bound the sampling cost of selecting a threshold from validation streams.

Admission commits a whole vehicle, whereas the plan assigns a fraction $y_i$ to the observed request. The next result connects the admission rule to the planned service fraction.

\begin{proposition}[Service threshold and spatial opportunity cost]
\label{prop:ts-admission}
For $\theta>0$, let $y_i^*$ denote the optimal planned service fraction of the current request in an exact solution of \eqref{eq:quadratic-plan}. Then \eqref{eq:price-rule} accepts if and only if $y_i^*\ge\gamma$, equivalently $g_i+\theta(1-2\gamma)\ge0$. Every set of optimal vehicle multipliers gives the same decision.
\end{proposition}

\begin{proof}
See Section~\ref{ec:price-rule} of the e-companion.\quad$\square$
\end{proof}

Proposition~\ref{prop:ts-admission} gives the admission rule a direct interpretation in terms of the service plan: the observed request is accepted when its planned fraction $y_i^*$ reaches the chosen threshold $\gamma$. For $\theta>0$, this regularized threshold admission rule links the unique planned fraction to a binary decision. Comparing forecast-induced changes in $y_i^*$ with its distance from $\gamma$ then allows us to quantify decision stability and revenue differences in Section~\ref{sec:theory}. Different optimal dual solutions may yield the same admission decision when $\theta>0$, even if their individual prices differ. Physically, accepting a ride exchanges a vehicle at the origin now for revenue and a vehicle at the destination later. A valuable destination therefore offsets part of the cost of leaving the origin. The threshold $\gamma$ determines how much planned service is required before the policy commits a whole vehicle. Lowering $\gamma$ accepts additional rides at the current prices.

The opportunity cost inferred from the mean-demand plan is only an approximation to the expected future revenue lost by serving the current request. The plan replaces random future requests with expected counts and allows fractional service, whereas actual admissions commit whole vehicles to realized requests. A zero adjusted-reward cutoff ($\gamma=1/2$) can therefore make suboptimal admission decisions. Calibrating $\gamma$ on independent validation streams selects the cutoff based on realized policy revenue. This adjusts the admission rule for these planning approximations and may improve revenue relative to the zero-cutoff rule.

\subsection{Fleet update and replanning}\label{sec:fleet-update}

If request $i$ is accepted, the policy assigns the idle vehicle at its origin with the smallest index $c$. The assigned vehicle's records for the next request $i+1$ become $s_{c,i+1}=d_{\xi_i}$ and $u_{c,i+1}=t_i+\ell_{\xi_i}\Delta$. The records of the other vehicles remain unchanged. If the request is rejected, the entire fleet state remains unchanged. In either case, the request is recorded for demand estimation before the policy considers request $i+1$ at its arrival time $t_{i+1}$. The interarrival time $t_{i+1}-t_i$ follows the request stream, independently of the planning interval length $\Delta$. Only the current admission $a_i$ is implemented. When the next feasible request arrives, the policy solves a new plan using the updated state.

\section{Calibration and revenue guarantees}\label{sec:theory}

For a given $\theta>0$, we study how threshold calibration improves the decisions made from a mean-demand fleet plan. Section~\ref{sec:calibration-value} discusses optimality on a general network and derives the optimal threshold interval on a single route. Section~\ref{sec:common-threshold} bounds the revenue loss of the best common threshold on a general network and evaluates that loss when several origins share return demand at a hub. Section~\ref{sec:calibration-selection} gives a revenue guarantee for selecting a threshold from validation streams in the full fleet model. Finally, Sections~\ref{sec:forecast-support} and~\ref{sec:forecast-rate} quantify how errors in estimating demand affect admission decisions and revenue at a fixed threshold.

Throughout, let $R:=\max_{z\in\mathcal Z}r_z$ be the maximum per-request reward. We write $\widehat\pi_{\theta,\gamma}$ for the LSOA policy using the estimated counts $\widehat q_{z,\tau|i}$, with admission threshold $\gamma$ and regularization parameter $\theta$. We write $\pi_{\theta,\gamma}^{\mathbb P}$ for the same planning and admission rule using the true conditional mean counts $q_{z,\tau|i}^{\mathbb P}$. These conditional means exclude the current request and condition on the history and arrivals observed through request $i$. Both policies use exact solves and the same vehicle-assignment conventions. The reference policy $\pi_{\theta,\gamma}^{\mathbb P}$ still uses mean-demand planning; it is not generally the optimal stochastic policy.

\subsection{Optimality analysis for a common threshold}\label{sec:calibration-value}

On a general network, a common threshold is optimal if and only if it makes an optimal decision at every feasible request the fleet encounters, with probability one. Whenever accepting is strictly better, the planned fraction $y_i^{\mathbb P,*}$ must be at least $\gamma$. Whenever rejecting is strictly better, it must be below $\gamma$. Either decision is optimal at a tie. These conditions alone do not give a formula for $\gamma$ in terms of demand and rewards. Each dispatch changes where and when vehicles become available, which affects both the planned fraction and the value of accepting the next request.

We can derive the optimal thresholds explicitly on a single route. The operator has one idle vehicle and receives a request with reward $r_{\mathrm{low}}$ at time zero. All later requests arrive on the same route during $(0,\Delta)$ according to a Poisson process with a deterministic, possibly time-varying rate. They pay $r_{\mathrm{high}}$ each and have expected count $q$. Each trip lasts $\Delta$, so the vehicle can serve only one request. The rolling plan covers the remaining horizon up to $H=2\Delta$ and uses the same threshold for every request.

\begin{proposition}[Optimal thresholds on a single route]\label{prop:calibration-value}
Suppose $0<q<1$, $r_{\mathrm{low}}>\theta>0$, and $r_{\mathrm{high}}-r_{\mathrm{low}}>2\theta q$. Let $w=r_{\mathrm{high}}(1-e^{-q})$ be the expected revenue from waiting for a later ride. The set of thresholds that attain optimal expected revenue is
\[
\begin{cases}
(0,\,1-q], & r_{\mathrm{low}}>w,\\
(1-q,\,1/(1+q)], & r_{\mathrm{low}}<w,\\
(0,\,1/(1+q)], & r_{\mathrm{low}}=w.
\end{cases}
\]
\end{proposition}

\begin{proof}
See Section~\ref{ec:calibration-value} of the e-companion, which also gives the conditions for a fleet of any size.\quad$\square$
\end{proof}

The two endpoints arise from different decisions. The initial plan reserves $q$ of the vehicle for later rides and assigns $1-q$ to the current ride. To turn down that ride, the threshold must exceed $1-q$. Once a higher-reward request arrives, the fleet replans. If the vehicle is still available and $m$ more requests are expected, the plan assigns $1/(1+m)$ to the request in hand. Since $m\le q$, a threshold at most $1/(1+q)$ accepts the first higher-reward request whenever it arrives. A higher threshold can turn down an early request even when it turns out to be the last one of the day.

Proposition~\ref{prop:calibration-value} shows that one threshold can yield optimal admission decisions on a single route. Across several origins, the plan may assign the same fraction to two requests even when the fleet should accept one and reject the other. Section~\ref{sec:common-threshold} examines how much revenue the fleet loses in such cases, even with the best common threshold.

\subsection{Revenue loss under the best common threshold}\label{sec:common-threshold}

We measure this loss as the difference between the expected revenue under optimal dispatch and the highest expected revenue attainable with a common threshold. We first give a lower bound for a general network and then calculate the loss exactly when several origins share return demand at a hub.

To bound this loss on a general network, consider the first request at which the threshold differs from an optimal policy on the same request stream. Up to that point, the two policies face the same fleet state. Suppose a disagreement occurs with probability at least $p$ for every common threshold, and at the first disagreement the chosen action costs at least $d$ in expected revenue even with optimal dispatch afterward. Then even the best common threshold loses at least $pd$ in expectation. In practice, quantifying this bound is difficult. The probability of reaching a given state depends on earlier dispatches, and valuing the two actions there requires evaluating later requests throughout the network. Section~\ref{ec:network-threshold-loss} of the e-companion gives a bound based on one feasible reference policy.

We can calculate the loss explicitly when several origins share return demand at a hub. Each of $C$ identical origins has one idle vehicle and an initial request to the hub paying $r_{\mathrm{low}}$. As in Section~\ref{sec:calibration-value}, later requests on the same route pay $r_{\mathrm{high}}$ and arrive during $(0,\Delta)$, independently across origins, at constant Poisson rates with mean count $q$ per origin. All trips last $\Delta$, and there are no vehicles initially at the hub. At time $\Delta$, exactly $k$ return requests arrive there, each paying $r_{\mathrm{ret}}$. There are no other requests. Only vehicles sent at time zero arrive in time for these return trips. Holding a vehicle at its origin earns $r_{\mathrm{high}}(1-e^{-q})$ in expectation, so sending it now costs $\delta=r_{\mathrm{high}}(1-e^{-q})-r_{\mathrm{low}}$ in expected outbound revenue. The plan covers the remaining horizon up to $2\Delta$.

\begin{proposition}[Revenue loss in a shared return market]\label{prop:hub-threshold}
Suppose $0<k<C$, $0<q<1/2$, $r_{\mathrm{low}},r_{\mathrm{ret}}>\theta>0$, $r_{\mathrm{high}}-r_{\mathrm{low}}>2\theta q$, and $0<\delta<r_{\mathrm{ret}}$. The optimal policy sends exactly $k$ vehicles to the hub. A common threshold sends all $C$ vehicles or none. Its smallest expected revenue loss is
\[
V^*(\mathbb P)-\sup_{0<\gamma<1}V_{\pi_{\theta,\gamma}^{\mathbb P}}(\mathbb P)
=\min\{(C-k)\delta,\;k(r_{\mathrm{ret}}-\delta)\}.
\]
\end{proposition}

\begin{proof}
See Section~\ref{ec:hub-threshold} of the e-companion.\quad$\square$
\end{proof}

Each of the first $k$ vehicles sent to the hub earns a return reward that covers the outbound revenue forgone at its origin. Sending another vehicle costs $\delta$ and earns no additional return reward. The plan assigns the same fraction $1-q$ to every initial ride, so adjusting the threshold cannot send just the $k$ vehicles needed. Sending all vehicles costs $(C-k)\delta$ in excess dispatches. Keeping them all at their origins gives up $k(r_{\mathrm{ret}}-\delta)$ in net return revenue. If $k/C$ approaches a fraction strictly between zero and one as the fleet grows, even the smaller of these two losses grows in proportion to fleet size. Section~\ref{ec:hub-threshold} of the e-companion extends the analysis to different demand rates across origins and uncertain return demand.

Proposition~\ref{prop:hub-threshold} measures how much revenue the fleet would still lose if the operator knew which common threshold earned the most. This loss arises because a single threshold cannot dispatch only the vehicles needed at the hub. In practice, the operator must also select a threshold from validation data. Section~\ref{sec:calibration-selection} bounds how much expected revenue this choice can cost relative to the best candidate threshold.

\subsection{Finite-sample guarantees for threshold calibration}\label{sec:calibration-selection}

Fix $\theta>0$, the historical data used to initialize forecasting, and a nonempty finite set of thresholds $\mathcal G\subset(0,1)$ before validation. Conditional on this initialization, assume the $D^{\mathrm{val}}\ge1$ validation days and a new operating day are independent draws from the same distribution. Replay each candidate on every validation day, resetting the fleet and estimator at the start of each replay. Let $\widehat\gamma$ maximize mean daily validation revenue over $\mathcal G$, with a fixed tie rule. Each replay follows the candidate policy over the entire day, so its revenue includes the later rides made possible or infeasible by each dispatch. Write $W_{\theta,\gamma}^{\mathrm{val},d}$ for the realized revenue on day $d$. The selection criterion is the sample mean $\widehat V_{\theta,\gamma}^{\mathrm{val}}=(D^{\mathrm{val}})^{-1} \sum_{d=1}^{D^{\mathrm{val}}}W_{\theta,\gamma}^{\mathrm{val},d}$. Sharing streams across thresholds pairs the revenue comparisons within each day. Independence across days gives the following guarantee.

\begin{theorem}[Revenue guarantee for threshold calibration]\label{thm:calibration-selection}
Under the sampling conditions above, for any $\delta\in(0,1)$, with probability at least $1-\delta$ over the validation streams,
\begin{equation}\label{eq:calibration-selection}
V_{\widehat\pi_{\theta,\widehat\gamma}}(\mathbb P) \ge \max_{\gamma\in\mathcal G}V_{\widehat\pi_{\theta,\gamma}}(\mathbb P) -\underbrace{2RCT\sqrt{\frac{\log(2|\mathcal G|/\delta)}{2D^{\mathrm{val}}}}}_{\text{sampling penalty}}.
\end{equation}
Here the selected threshold is held fixed when evaluating revenue on the new operating day.
\end{theorem}

\begin{proof}
Conditional on the initialization, each candidate's daily revenue lies in $[0,RCT]$ because each vehicle can complete at most $T$ rides. Hoeffding's inequality and a union bound over $\mathcal G$ bound the largest absolute difference between mean validation revenue and expected revenue by half the penalty in~\eqref{eq:calibration-selection} with probability at least $1-\delta$. Comparing the empirical maximizer with an expected-revenue maximizer costs at most twice this difference.\quad$\square$
\end{proof}

For a fixed candidate set, the sampling penalty decreases at rate $(D^{\mathrm{val}})^{-1/2}$. If $1/2\in\mathcal G$ and the best candidate improves on the zero-cutoff rule by more than this penalty, then on the stated event the selected threshold also improves on that rule. The same guarantee applies to plans that use true conditional means. More validation days help identify the best cutoff in the set. The loss from using one cutoff across origins, quantified in Proposition~\ref{prop:hub-threshold}, can persist even when that cutoff is known.

\subsection{Forecast error after calibration}\label{sec:forecast-support}

Fix $\theta$ and $\gamma$. We compare $\widehat\pi_{\theta,\gamma}$ with $\pi_{\theta,\gamma}^{\mathbb P}$ on the same request stream. Both policies use the same fleet plan and threshold, so their revenue difference measures the effect of estimating demand.

At the same fleet state and planning grid, let $\widehat y_i^*$ and $y_i^{\mathbb P,*}$ be the planned current fractions under estimated and true conditional forecasts. Write $\widehat{\boldsymbol q}$ and $\boldsymbol q^{\mathbb P}$ for the corresponding forecast arrays, and let $\mathcal Z^{\mathrm{end}}_{\tau|i} :=\{z\in\mathcal Z:t_i+(\tau+\ell_z)\Delta\le H\}$ be the eligible request types that can finish by $H$. Their total forecast difference is defined as  $\varepsilon_i(\widehat{\boldsymbol q},\boldsymbol q^{\mathbb P}) :=\sum_{\tau\in\mathcal T_i}\sum_{z\in\mathcal Z^{\mathrm{end}}_{\tau|i}} |\widehat q_{z,\tau|i}-q_{z,\tau|i}^{\mathbb P}|$.

\begin{lemma}[Sensitivity to demand-estimation error]
\label{lem:calibration-forecast-support}
Suppose both plans use the same $\theta>0$ and finite nonnegative forecast counts. Then
\begin{equation}\label{eq:calibration-forecast-support}
|\widehat y_i^*-y_i^{\mathbb P,*}| \le\min\{1,\varepsilon_i(\widehat{\boldsymbol q},\boldsymbol q^{\mathbb P})\}.
\end{equation}
\end{lemma}

\begin{proof}
See Section~\ref{ec:sensitivity-proof} of the e-companion.\quad$\square$
\end{proof}

The bound is nontrivial when the total count error is less than one. The proof uses the vehicle-flow structure and normalized quadratic regularization to translate forecast changes into a linear bound on the current service fraction.

If $\varepsilon_i(\widehat{\boldsymbol q},\boldsymbol q^{\mathbb P}) <|y_i^{\mathbb P,*}-\gamma|$, the forecast error cannot move the planned fraction across $\gamma$, so the admission decision is unchanged. To extend this local observation to daily revenue, we evaluate these quantities along the reference trajectory of $\pi_{\theta,\gamma}^{\mathbb P}$. Indicators involving $|y_i^{\mathbb P,*}-\gamma|$ are defined to be zero for requests that are infeasible on that trajectory. The next result, Theorem~\ref{thm:ts-oracle}, translates these possible admission differences into a bound on the expected revenue difference between the two policies on a common request stream.

\begin{theorem}[Revenue comparison on a common request stream]
\label{thm:ts-oracle}
If $\mathbb E_{\mathbb P}[I^2]<\infty$, then
\begin{equation}\label{eq:ts-oracle-comparison}
|V_{\widehat\pi_{\theta,\gamma}}(\mathbb P)-V_{\pi_{\theta,\gamma}^{\mathbb P}}(\mathbb P)| \le R\,\mathbb E_{\mathbb P}\!\left[ \sum_{i=1}^{I}(I-i+1)\, \mathbbm{1}\left\{|y_i^{\mathbb P,*}-\gamma|\le \min\left\{1,\varepsilon_i(\widehat{\boldsymbol q},\boldsymbol q^{\mathbb P})\right\} \right\}\right].
\end{equation}
\end{theorem}

\begin{proof}
See Section~\ref{ec:revenue} of the e-companion.\quad$\square$
\end{proof}

The indicator in \eqref{eq:ts-oracle-comparison} marks requests whose reference fraction $y_i^{\mathbb P,*}$ is sufficiently close to $\gamma$ that forecast errors could change acceptance or rejection. The factor $I-i+1$ accounts for all requests remaining from that decision onward. The argument couples the policies on the same stream: before their first different decision, their fleet states and revenues agree; afterward, their remaining revenues differ by at most $R(I-i+1)$.

Thus forecast accuracy matters most near the admission threshold. Away from it, changing the forecast need not change the action. Theorem~\ref{thm:ts-oracle} imposes no margin condition or Poisson assumption, but its bound still depends on how often such near-threshold decisions occur. We next control this frequency to obtain a revenue convergence rate as historical data accumulate.

\subsection{Revenue convergence with historical data}\label{sec:forecast-rate}

To quantify how the revenue difference changes with the number of historical days, we specify the arrival process used for estimation.

\begin{assumption}[Poisson arrivals within demand contexts]
\label{ass:ts-poisson}
The request types and demand-context partition are fixed. We assume that arrivals of each type $z$ follow a Poisson process with deterministic intensity $\alpha_{z,\kappa(u)}$ at time $u$, where $0\le\alpha_{z,j}<\infty$ for every context $j$. Processes for different types are independent. Historical streams and the operating-day stream are independent and use the same intensities. Historical observation durations are fixed in advance.
\end{assumption}

Context rates may vary within a day but are stable across days. We also control how often the service fraction computed with these rates falls near the threshold.
\begin{assumption}[Limited concentration near the admission threshold]
\label{ass:revenue-margin}
For the reference policy $\pi_{\theta,\gamma}^{\mathbb P}$, we assume that there exists a constant $0<M_{\theta,\gamma}<\infty$ such that, for every $n\ge1$ with $\mathbb P(I=n)>0$ and every $m\in[0,1]$,
\begin{equation}\label{eq:weighted-margin}
\mathbb E_{\mathbb P}\!\left[ \sum_{i=1}^{n}(n-i+1)\mathbbm{1}\{|y_i^{\mathbb P,*}-\gamma|\le m\}\,\middle|\,I=n \right] \le M_{\theta,\gamma}\,m\,\frac{n(n+1)}{2}.
\end{equation}
The constant is independent of the number of historical days $D$ and the daily request count $n$.
\end{assumption}

This assumption requires that the reference policy's planned service fractions do not concentrate too heavily near the admission threshold $\gamma$. Specifically, conditional on $I=n$, the expected weighted share of reference decisions within distance $m$ of $\gamma$ (the quantity $\mathbb E_{\mathbb P}[\sum_{i=1}^{n} \frac{2(n-i+1)}{n(n+1)}\mathbbm{1}\{|y_i^{\mathbb P,*}-\gamma|\le m\} \mid I=n]$) is at most $M_{\theta,\gamma}m$. The weights are normalized by their total, $\sum_{i=1}^{n}(n-i+1)=n(n+1)/2$.

This assumption may be reasonable when requests differ substantially in revenue, trip duration, and destination demand, so that most admission decisions are clear and only a small share are close to the threshold $\gamma$. Changes in demand and vehicle availability may also spread the reference fractions $y_i^{\mathbb P,*}$ across different values. These features can support the assumption, but do not guarantee it.

Let $h_j:=\int_0^H\mathbbm{1}\{\kappa(u)=j\}\,du$ be the daily duration of context $j$. Set $h_{\min}:=\min_jh_j$ and let $\Lambda:=\sum_{z\in\mathcal Z}\sum_{j\in\mathcal J}\alpha_{z,j}h_j$ be the expected daily request count. The next result combines forecast accuracy with Assumption~\ref{ass:revenue-margin} to bound the revenue cost of estimation.

\begin{theorem}[Revenue convergence with historical data]
\label{thm:revenue-rate}
Under Assumptions~\ref{ass:ts-poisson} and~\ref{ass:revenue-margin}, suppose each of $D\ge1$ historical days is observed throughout $[0,H)$ with the same context mapping $\kappa(u)$ as the operating day, with $h_j>0$ and $E_j^{\mathrm{hist}}=Dh_j$ for every $j\in\mathcal J$. With a fixed initial fleet and fixed parameters $\theta>0$ and $\gamma\in(0,1)$,
\begin{equation}\label{eq:revenue-learning-rate}
|V_{\widehat\pi_{\theta,\gamma}}(\mathbb P) -V_{\pi_{\theta,\gamma}^{\mathbb P}}(\mathbb P)| \le R M_{\theta,\gamma}L\Delta \left[ \frac{\Lambda^2+2\Lambda}{2\sqrt D} \sum_{z\in\mathcal Z}\sum_{j\in\mathcal J} \sqrt{\frac{\alpha_{z,j}}{h_j}} +\frac{\Lambda^3+3\Lambda^2+\Lambda}{D h_{\min}} \right] =O(D^{-1/2}),
\end{equation}
with all model parameters held fixed as $D\to\infty$.
\end{theorem}

\begin{proof}
See Section~\ref{ec:revenue-rate} of the e-companion.\quad$\square$
\end{proof}

Unlike \eqref{eq:ts-oracle-comparison}, which leaves the frequency of near-threshold decisions unspecified, Assumption~\ref{ass:revenue-margin} controls that frequency and yields an explicit revenue rate.

The bound shows how sample size and decision sensitivity interact: more historical days reduce the loss, while a larger $M_{\theta,\gamma}$ allows more reference decisions near $\gamma$ and weakens the guarantee. The rate $D^{-1/2}$ is tight in general under these assumptions; Section~\ref{ec:revenue-rate-tightness} gives a fixed instance with a $\Theta(D^{-1/2})$ revenue difference. Without Assumption~\ref{ass:revenue-margin}, Theorem~\ref{thm:ts-oracle} remains valid, but small forecast errors alone do not guarantee the revenue rate in Theorem~\ref{thm:revenue-rate}. In particular, positive probability at $y_i^{\mathbb P,*}=\gamma$ can leave admissions sensitive to arbitrarily small estimation errors.

% Source: frozen Gurobi study.
\section{Case study}\label{sec:case-study}

We compare LSOA with seven benchmarks across six synthetic demand families and a real-world case study based on New York City trip records.

\subsection{Case study design}\label{sec:case-design}

\subsubsection{Synthetic demand.}
\begin{figure}[htp]
\centering
\includegraphics[width=\textwidth]{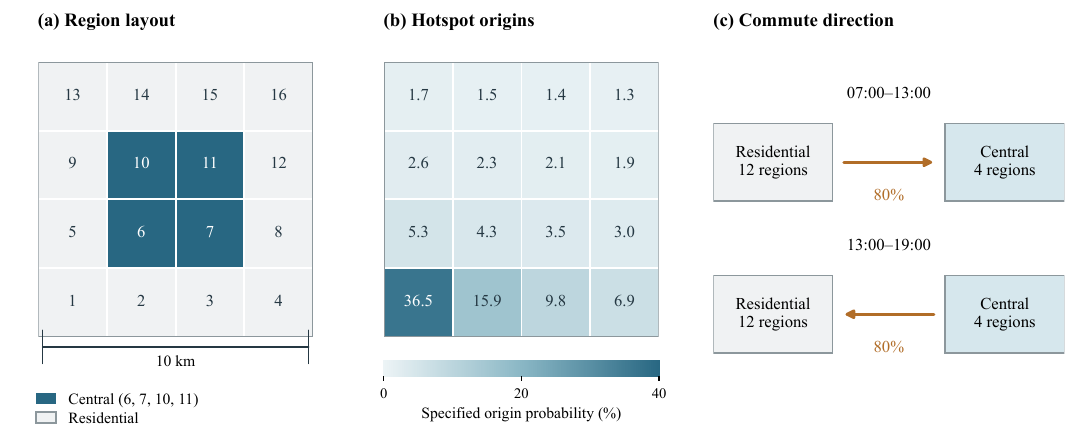}
\caption{Synthetic service area and spatial demand.}\label{fig:synthetic-map}
\end{figure}

Figure~\ref{fig:synthetic-map} shows the synthetic service area, a $10\times10$ km square divided into $S=16$ equal regions with $C=100$ vehicles. Operations run from 07:00 to 19:00. A request type $z$ is defined by its ordered origin--destination pair $(o_z,d_z)$, including trips within the same region, giving $Z=16^2=256$ types. Requests of the same type share reward $r_z$ and duration $\ell_z\Delta$. Trip distances use Manhattan distance between region centers, with a minimum of 1 km. Travel speed is 20 km/h, durations are rounded up to two-minute intervals, and rewards are \$3 plus \$2 per kilometer.

\begin{figure}[htp]
\centering
\includegraphics[width=\textwidth]{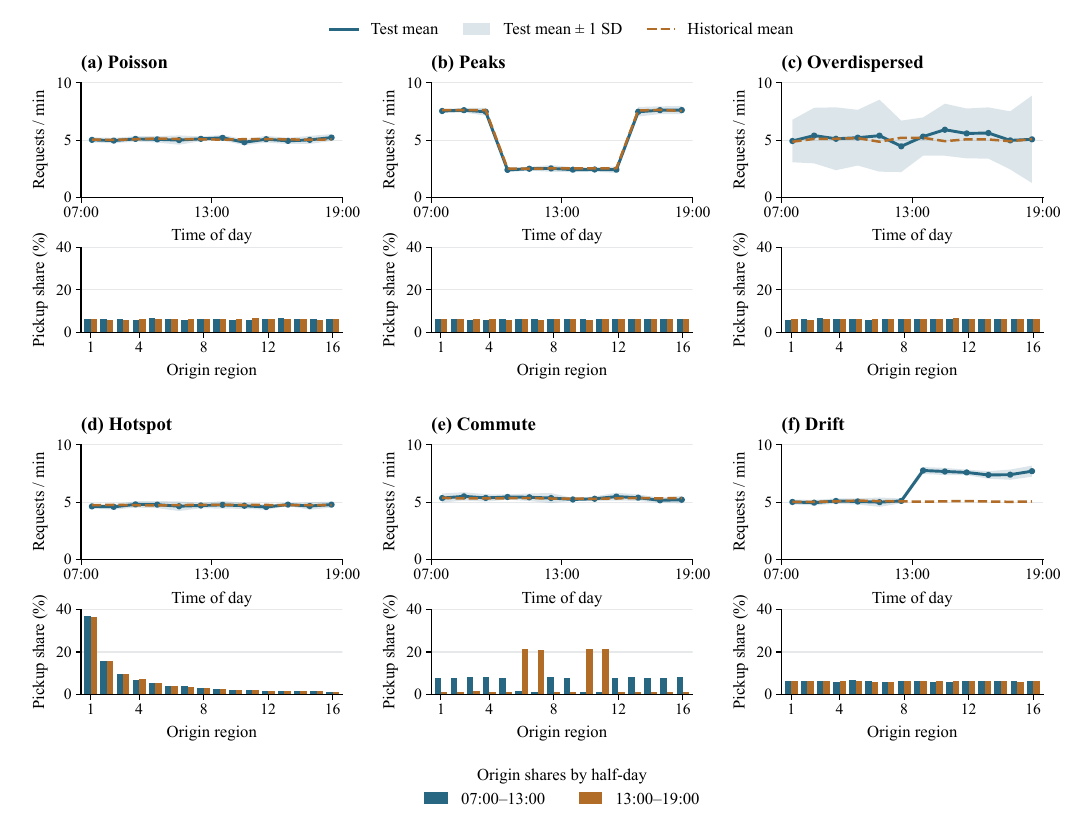}
\caption{Synthetic demand patterns.}\label{fig:synthetic-demand}
\end{figure}

The six demand families in Figure~\ref{fig:synthetic-demand} vary temporal demand and spatial flows. Poisson uses a constant arrival rate. Peaks uses 1.5 times the daily mean rate during 07:00--10:00 and 16:00--19:00, and 0.5 times otherwise. Overdispersed uses negative binomial hourly counts with variance equal to the mean plus the squared mean divided by five. Drift raises the test-day rate by 50\% after 13:00 relative to history. These four families use uniform origins and destinations. Hotspot assigns origin probabilities proportional to $s^{-1.2}$, where $s$ is the region ID, with uniformly distributed destinations. Commute sends 80\% of trips from the twelve residential regions to the four central regions (6, 7, 10, and 11) before 13:00 and reverses the direction afterward. Its remaining 20\% have uniform origins and destinations. Figure~\ref{fig:synthetic-demand} reports mean hourly arrival rates and regional pickup shares, with shading showing one sample standard deviation (SD) of test-day rates.

For each family, we choose the base arrival rate so that the expected number of arrivals per hour, multiplied by the mean trip duration in hours, equals $C$. We use ten instances per family, each with 20 historical days and one test day. Threshold calibration uses two separate validation days per family in the main comparison; we retain 20 validation days for the sample-size analysis. Test days contain 3,301--4,688 requests. All policies share the same history, test requests, and initially idle fleet, allocated by historical origin shares. Type-level rewards and durations remain fixed across days.

Threshold calibration maximizes average revenue on two validation streams with the first instance's historical sample. Ties are resolved in favor of the larger threshold. The selected threshold is then held fixed across that family's ten test instances, each retaining its original 20-day historical sample. Validation streams are independent of the test streams, and each replay resets the fleet and demand estimator. Drift uses pre-drift validation streams and retains the test-day demand increase.

\subsubsection{New York requests.}
\begin{figure}[htp]
\centering
\includegraphics[width=\textwidth]{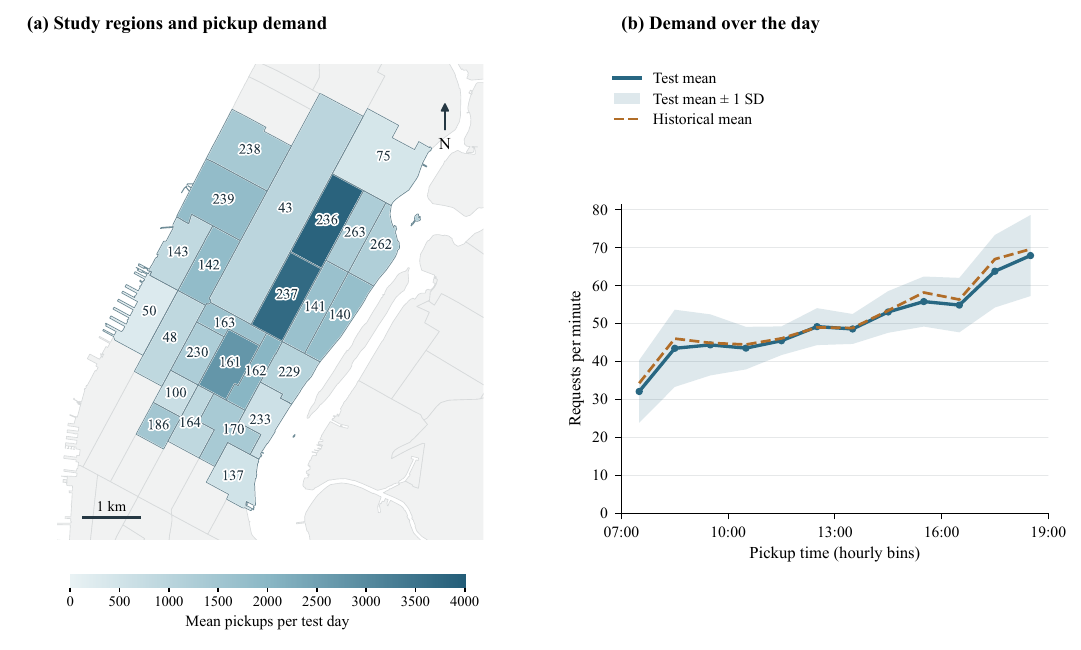}
\caption{New York study regions and demand.}\label{fig:nyc-demand}
\end{figure}

Figure~\ref{fig:nyc-demand} shows the New York study area and demand patterns. The area comprises $S=25$ adjacent Manhattan taxi zones around Midtown Center, listed in Section~\ref{ec:nyc-service-area} of the e-companion. As in the synthetic case, each ordered pair of pickup and drop-off zones $(o_z,d_z)$ defines a request type $z$. Including within-zone trips gives $Z=25^2=625$ types. We use January and February 2025 yellow-taxi records from the New York City Taxi and Limousine Commission \citep{nyctlc2025data}. Each pickup is treated as a request at its recorded time. We include trips that start and end within the study area and depart between 07:00 and 19:00. We retain trips with recorded durations of 1--120 minutes and fares greater than zero and no more than \$200. Trips with the same pickup time are processed in the order recorded in the data.

The historical sample contains 370,872 requests from the ten weekdays of January 6--17. January 21--22 are used for validation, and the ten test dates are February 6--7, 10--14, and 17--19. Test days contain 23,970--39,121 requests, averaging 50.18 requests per minute over the 12-hour operating day. Figure~\ref{fig:nyc-demand}(a) shows mean daily pickups by zone, and panel (b) compares mean hourly arrival rates across the ten historical and ten test days. Shading shows one sample standard deviation across test days.

All requests with the same origin--destination pair share reward $r_z$ and travel time $\ell_z\Delta$, determined by that pair's historical median fare and duration, with durations rounded up to two-minute intervals. Rewards range from \$5.80 to \$31.16 and durations from 4 to 32 minutes. Revenue sums these fixed type-level rewards across accepted requests. We test fleets of $C=250$, $500$, and $750$ vehicles, initialized according to historical origin shares. We select $\gamma$ separately for each fleet size using January 21 and 22, over the same 19-point candidate set $\mathcal G=\{0.05,0.10,\ldots,0.95\}$ as in the synthetic study. Ties favor the larger threshold. Validation days are excluded from the forecasting history.

\subsection{Benchmark methods}\label{sec:benchmarks}
The benchmarks use the same request stream and initial fleet as LSOA. Additional implementation details are in Section~\ref{ec:benchmark-models}.

\noindent\textbf{Greedy admission.} For each physically feasible request $i$, Greedy uses $a_i=\mathbbm{1}\{r_{\xi_i}\ge0\}$, without reserving capacity for future demand. Since rewards are nonnegative in both studies, it accepts every physically feasible request.

\noindent\textbf{Mean Demand Lookahead.} This policy solves \eqref{eq:quadratic-plan} with $\theta=0$ twice, fixing $y_i=1$ for acceptance and $y_i=0$ for rejection. Both solves use the same demand forecasts and 30-minute planning grid as LSOA. Let $W_i^{\mathrm{LP}}(a)$ be the optimized future revenue under action $a$, excluding the current fare. The policy accepts when $r_{\xi_i}\ge W_i^{\mathrm{LP}}(0)-W_i^{\mathrm{LP}}(1)$.

\noindent\textbf{Sampled Lookahead.} At each feasible arrival, this policy generates 20 future request sequences for the next 30 minutes. It samples arrival times for each origin--destination type from Poisson processes with learned rates $\widehat\alpha_{z,j|i}$; each sampled request inherits its type's reward $r_z$ and duration $\ell_z\Delta$. For each sequence, it adapts \eqref{eq:quadratic-plan} with $\theta=0$, replacing interval-level expected demand by individual sampled requests with binary service decisions at the sampled arrival and completion times. It solves this plan with $y_i=1$ and $y_i=0$ to obtain the optimized future revenues after acceptance and rejection, excluding the current fare. It uses the same admission test as Mean Demand Lookahead, replacing the value difference by its average over the 20 paired sequences. Each optimization knows its entire sampled sequence, so this average need not recover the optimal online value function.

\noindent\textbf{Supply-aware Mirror Backpressure (MBP).} Following \citet{kanoria2022mirror}, this policy combines regional idle-vehicle counts with a supply price for vehicle occupation. For a physically feasible request, it accepts when $r_{\xi_i}+f(\bar n_{o_{\xi_i},i})-f(\bar n_{d_{\xi_i},i}) -v_i^-\ell_{\xi_i}\Delta\ge0$. Here, $\bar n_{s,i}$ is a normalized idle-vehicle count in region $s$, $f$ is an increasing congestion function, and $v_i^-$ is the supply price per vehicle-minute immediately before admission. The supply price increases with accepted trip duration and decreases with elapsed time. This policy requires neither demand forecasts nor optimization. Detailed definitions of $\bar n_{s,i}$, $f$, and $v_i^-$ are given in Section~\ref{ec:supply-aware-mbp}.

\noindent\textbf{Batch.} Every two minutes, let $\mathcal E$ contain vehicle--request pairs $(c,i)$ with an idle vehicle at the waiting request's origin and completion by $H$. Batch maximizes $\sum_{(c,i)\in\mathcal E}r_{\xi_i}a_{c,i}^{\mathrm{Batch}}$ over binary assignments $a_{c,i}^{\mathrm{Batch}}\in\{0,1\}$, subject to $\sum_{i:(c,i)\in\mathcal E}a_{c,i}^{\mathrm{Batch}}\le1$ for each vehicle and $\sum_{c:(c,i)\in\mathcal E}a_{c,i}^{\mathrm{Batch}}\le1$ for each waiting request. Unmatched requests leave at dispatch. Accepted rides start then, shifting completion times by the waiting time; passenger waits are reported alongside revenue.

\noindent\textbf{Upper bound.} Knowing the full realized request stream, the offline model computes $\max_{\{a_i^{\mathrm{off}}\}}\allowbreak\sum_{i\in\mathcal I}r_{\xi_i}a_i^{\mathrm{off}}$, where $a_i^{\mathrm{off}}\in\{0,1\}$ indicates service at the fixed arrival time $t_i$. It uses the vehicle-flow accounting in \eqref{eq:quadratic-vehicle} at realized arrival and completion times, with $a_i^{\mathrm{off}}=0$ if $t_i+\ell_{\xi_i}\Delta>H$. It bounds revenue for immediate-start policies.

\subsection{Model configurations and evaluation protocol}\label{sec:model-configurations}
All LSOA configurations solve \eqref{eq:quadratic-plan} to get dual prices and use \eqref{eq:price-rule} for admission decisions. The primary configuration uses positive regularization $\theta>0$ and a validation-selected admission threshold $\gamma$.

LSOA ($\theta=0$) is listed as a separate method in the result tables. It uses the same forecasts, planning horizon, and vehicle constraints, but solves the planning problem without regularization. Its admission rule is \(a_i=\mathbbm{1}\{g_i\ge0\}\), with \(g_i\) calculated using the dual prices from this linear problem.

We process requests in arrival order. Sequential policies accept a request only if an idle vehicle is available in its origin region and the trip can finish by 19:00. Accepted rides start immediately and release their vehicles after the calibrated travel time. The simulation excludes empty relocation and travel within a pickup region. Batch uses the same fleet but allows requests to wait until the next two-minute dispatch epoch.

LSOA, Mean Demand Lookahead, and Sampled Lookahead use a 30-minute planning window. LSOA and Mean Demand Lookahead use up to $L=15$ planning intervals of length $\Delta=2$ minutes starting at request time $t_i$. In both the synthetic and New York studies, we divide the 12-hour operating day into one-hour demand contexts, giving $J=12$ contexts in total. We index them in chronological order: $j=1$ corresponds to 07:00--08:00, $j=2$ to 08:00--09:00, and so on through $j=12$ for 18:00--19:00. Each context $j$ has a separate arrival-rate estimate $\widehat\alpha_{z,j|i}$ for each origin--destination type $z$, using observations from the same hour across days. These hourly contexts capture time-of-day demand patterns, while the two-minute planning intervals determine when future vehicle movements are scheduled. We update demand rates $\widehat\alpha_{z,j|i}$ using \eqref{eq:demand-estimate}, adding every observed request after its decision, including rejected and infeasible requests. For LSOA and Mean Demand Lookahead, the planning grid and vehicle releases follow Section~\ref{sec:planning}.

All optimization problems are solved with Gurobi 13.0.3 on a machine with an AMD Ryzen 9 9950X3D processor (16 physical cores, 32 threads) and 64~GB of installed memory. The implementation uses the sparse flow model and combines parallel arcs as described in Sections~\ref{ec:sparse-primal} and~\ref{ec:parallel-arcs}. Solver settings and verification checks are given in Section~\ref{ec:implementation}.

\subsection{Synthetic results}\label{sec:case-results}

For the main synthetic comparison, we fix $\theta=10$ and $D^{\mathrm{val}}=2$. We calibrate $\gamma$ separately for each demand family over $\mathcal G=\{0.05,0.10,\ldots,0.95\}$. We also report LSOA ($\gamma=1/2$) with the same $\theta=10$; this is the zero-cutoff rule $g_i\ge0$ and isolates the effect of threshold calibration.

\subsubsection{Revenue and computation time.}
Table~\ref{tab:synthetic-performance} reports daily revenue across the six demand families. Two validation days select $\gamma=0.15$ for Peaks, $\gamma=0.45$ for Commute, and $\gamma=0.25$ for the other four families.

\begin{table}[!ht]
\centering
\caption{Synthetic revenue.}\label{tab:synthetic-performance}
\fontsize{7.5}{9}\selectfont
\renewcommand{\scriptsize}{\fontsize{6.8}{8}\selectfont}
\setlength{\tabcolsep}{1pt}
\renewcommand{\arraystretch}{1.15}
\begin{tabular}{@{}lccccccccc@{}}
\hline\noalign{\vskip 4pt}
Demand & Greedy & LSOA & \shortstack{LSOA\\($\gamma=1/2$)} & \shortstack{LSOA\\($\theta=0$)} & \shortstack{Mean Demand\\Lookahead} & \shortstack{Sampled\\Lookahead} & \shortstack{Supply-aware\\MBP} & Batch & Upper bound \\
\hline\noalign{\vskip 4pt}
Poisson & $37.75 \pm 0.49$ & \shortstack[c]{$38.59 \pm 0.56$ \\ {\scriptsize \textcolor[rgb]{0.08,0.42,0.23}{$\uparrow 2.21 \pm 0.64\%$}}} & \shortstack[c]{$36.69 \pm 0.46$ \\ {\scriptsize \textcolor[rgb]{0.70,0.15,0.13}{$\downarrow 2.81 \pm 0.54\%$}}} & \shortstack[c]{$38.60 \pm 0.63$ \\ {\scriptsize \textcolor[rgb]{0.08,0.42,0.23}{$\uparrow 2.25 \pm 0.80\%$}}} & \shortstack[c]{$36.53 \pm 0.55$ \\ {\scriptsize \textcolor[rgb]{0.70,0.15,0.13}{$\downarrow 3.24 \pm 0.84\%$}}} & \shortstack[c]{$38.35 \pm 0.54$ \\ {\scriptsize \textcolor[rgb]{0.08,0.42,0.23}{$\uparrow 1.58 \pm 0.59\%$}}} & \shortstack[c]{$37.26 \pm 0.41$ \\ {\scriptsize \textcolor[rgb]{0.70,0.15,0.13}{$\downarrow 1.31 \pm 0.67\%$}}} & \shortstack[c]{$38.71 \pm 0.44$ \\ {\scriptsize \textcolor[rgb]{0.08,0.42,0.23}{$\uparrow 2.54 \pm 0.50\%$}}} & $42.77 \pm 0.49$ \\[1pt]
Peaks & $33.56 \pm 0.58$ & \shortstack[c]{$33.79 \pm 0.66$ \\ {\scriptsize \textcolor[rgb]{0.08,0.42,0.23}{$\uparrow 0.69 \pm 0.56\%$}}} & \shortstack[c]{$31.32 \pm 0.52$ \\ {\scriptsize \textcolor[rgb]{0.70,0.15,0.13}{$\downarrow 6.67 \pm 0.94\%$}}} & \shortstack[c]{$33.58 \pm 0.58$ \\ {\scriptsize \textcolor[rgb]{0.08,0.42,0.23}{$\uparrow 0.05 \pm 0.86\%$}}} & \shortstack[c]{$31.33 \pm 0.55$ \\ {\scriptsize \textcolor[rgb]{0.70,0.15,0.13}{$\downarrow 6.66 \pm 0.73\%$}}} & \shortstack[c]{$33.28 \pm 0.61$ \\ {\scriptsize \textcolor[rgb]{0.70,0.15,0.13}{$\downarrow 0.84 \pm 0.64\%$}}} & \shortstack[c]{$32.99 \pm 0.51$ \\ {\scriptsize \textcolor[rgb]{0.70,0.15,0.13}{$\downarrow 1.70 \pm 0.70\%$}}} & \shortstack[c]{$34.55 \pm 0.63$ \\ {\scriptsize \textcolor[rgb]{0.08,0.42,0.23}{$\uparrow 2.95 \pm 0.62\%$}}} & $37.06 \pm 0.56$ \\[1pt]
Overdispersed & $36.02 \pm 1.67$ & \shortstack[c]{$36.59 \pm 1.62$ \\ {\scriptsize \textcolor[rgb]{0.08,0.42,0.23}{$\uparrow 1.59 \pm 0.71\%$}}} & \shortstack[c]{$34.80 \pm 1.49$ \\ {\scriptsize \textcolor[rgb]{0.70,0.15,0.13}{$\downarrow 3.35 \pm 1.27\%$}}} & \shortstack[c]{$36.57 \pm 1.65$ \\ {\scriptsize \textcolor[rgb]{0.08,0.42,0.23}{$\uparrow 1.53 \pm 0.64\%$}}} & \shortstack[c]{$34.58 \pm 1.58$ \\ {\scriptsize \textcolor[rgb]{0.70,0.15,0.13}{$\downarrow 3.99 \pm 1.29\%$}}} & \shortstack[c]{$36.35 \pm 1.61$ \\ {\scriptsize \textcolor[rgb]{0.08,0.42,0.23}{$\uparrow 0.93 \pm 0.85\%$}}} & \shortstack[c]{$35.33 \pm 1.47$ \\ {\scriptsize \textcolor[rgb]{0.70,0.15,0.13}{$\downarrow 1.90 \pm 0.78\%$}}} & \shortstack[c]{$36.97 \pm 1.65$ \\ {\scriptsize \textcolor[rgb]{0.08,0.42,0.23}{$\uparrow 2.67 \pm 0.89\%$}}} & $40.36 \pm 1.81$ \\[1pt]
Hotspot & $15.45 \pm 0.57$ & \shortstack[c]{$17.15 \pm 0.83$ \\ {\scriptsize \textcolor[rgb]{0.08,0.42,0.23}{$\uparrow 10.99 \pm 2.96\%$}}} & \shortstack[c]{$16.34 \pm 0.72$ \\ {\scriptsize \textcolor[rgb]{0.08,0.42,0.23}{$\uparrow 5.79 \pm 4.16\%$}}} & \shortstack[c]{$17.35 \pm 0.67$ \\ {\scriptsize \textcolor[rgb]{0.08,0.42,0.23}{$\uparrow 12.29 \pm 2.61\%$}}} & \shortstack[c]{$16.83 \pm 0.65$ \\ {\scriptsize \textcolor[rgb]{0.08,0.42,0.23}{$\uparrow 8.96 \pm 3.28\%$}}} & \shortstack[c]{$16.93 \pm 0.75$ \\ {\scriptsize \textcolor[rgb]{0.08,0.42,0.23}{$\uparrow 9.59 \pm 3.37\%$}}} & \shortstack[c]{$15.37 \pm 0.51$ \\ {\scriptsize \textcolor[rgb]{0.70,0.15,0.13}{$\downarrow 0.54 \pm 0.92\%$}}} & \shortstack[c]{$14.93 \pm 0.75$ \\ {\scriptsize \textcolor[rgb]{0.70,0.15,0.13}{$\downarrow 3.39 \pm 1.97\%$}}} & $23.62 \pm 0.73$ \\[1pt]
Commute & $9.43 \pm 0.33$ & \shortstack[c]{$12.31 \pm 0.53$ \\ {\scriptsize \textcolor[rgb]{0.08,0.42,0.23}{$\uparrow 30.61 \pm 3.62\%$}}} & \shortstack[c]{$12.24 \pm 0.52$ \\ {\scriptsize \textcolor[rgb]{0.08,0.42,0.23}{$\uparrow 29.77 \pm 3.56\%$}}} & \shortstack[c]{$11.02 \pm 0.48$ \\ {\scriptsize \textcolor[rgb]{0.08,0.42,0.23}{$\uparrow 16.88 \pm 2.00\%$}}} & \shortstack[c]{$12.05 \pm 0.46$ \\ {\scriptsize \textcolor[rgb]{0.08,0.42,0.23}{$\uparrow 27.77 \pm 2.64\%$}}} & \shortstack[c]{$12.36 \pm 0.51$ \\ {\scriptsize \textcolor[rgb]{0.08,0.42,0.23}{$\uparrow 31.11 \pm 3.76\%$}}} & \shortstack[c]{$9.41 \pm 0.33$ \\ {\scriptsize \textcolor[rgb]{0.70,0.15,0.13}{$\downarrow 0.21 \pm 0.51\%$}}} & \shortstack[c]{$9.64 \pm 0.39$ \\ {\scriptsize \textcolor[rgb]{0.08,0.42,0.23}{$\uparrow 2.21 \pm 1.14\%$}}} & $14.41 \pm 0.49$ \\[1pt]
Drift & $41.43 \pm 0.42$ & \shortstack[c]{$41.98 \pm 0.41$ \\ {\scriptsize \textcolor[rgb]{0.08,0.42,0.23}{$\uparrow 1.31 \pm 0.41\%$}}} & \shortstack[c]{$39.69 \pm 0.45$ \\ {\scriptsize \textcolor[rgb]{0.70,0.15,0.13}{$\downarrow 4.21 \pm 0.56\%$}}} & \shortstack[c]{$42.05 \pm 0.45$ \\ {\scriptsize \textcolor[rgb]{0.08,0.42,0.23}{$\uparrow 1.48 \pm 0.46\%$}}} & \shortstack[c]{$39.59 \pm 0.34$ \\ {\scriptsize \textcolor[rgb]{0.70,0.15,0.13}{$\downarrow 4.45 \pm 0.48\%$}}} & \shortstack[c]{$41.64 \pm 0.39$ \\ {\scriptsize \textcolor[rgb]{0.08,0.42,0.23}{$\uparrow 0.51 \pm 0.38\%$}}} & \shortstack[c]{$40.63 \pm 0.36$ \\ {\scriptsize \textcolor[rgb]{0.70,0.15,0.13}{$\downarrow 1.94 \pm 0.70\%$}}} & \shortstack[c]{$42.71 \pm 0.46$ \\ {\scriptsize \textcolor[rgb]{0.08,0.42,0.23}{$\uparrow 3.09 \pm 0.45\%$}}} & $46.06 \pm 0.37$ \\[1pt]
\hline\noalign{\vskip 4pt}
Average & $28.94 \pm 12.15$ & \shortstack[c]{$\boldsymbol{30.07 \pm 11.33}$ \\ {\scriptsize \textcolor[rgb]{0.08,0.42,0.23}{$\boldsymbol{\uparrow 7.90 \pm 11.00\%}$}}} & \shortstack[c]{$28.51 \pm 10.54$ \\ {\scriptsize \textcolor[rgb]{0.08,0.42,0.23}{$\uparrow 3.09 \pm 12.85\%$}}} & \shortstack[c]{$29.86 \pm 11.64$ \\ {\scriptsize \textcolor[rgb]{0.08,0.42,0.23}{$\uparrow 5.75 \pm 6.62\%$}}} & \shortstack[c]{$28.48 \pm 10.44$ \\ {\scriptsize \textcolor[rgb]{0.08,0.42,0.23}{$\uparrow 3.07 \pm 12.38\%$}}} & \shortstack[c]{$29.82 \pm 11.21$ \\ {\scriptsize \textcolor[rgb]{0.08,0.42,0.23}{$\uparrow 7.15 \pm 11.51\%$}}} & \shortstack[c]{$28.50 \pm 11.86$ \\ {\scriptsize \textcolor[rgb]{0.70,0.15,0.13}{$\downarrow 1.27 \pm 0.97\%$}}} & \shortstack[c]{$29.59 \pm 12.70$ \\ {\scriptsize \textcolor[rgb]{0.08,0.42,0.23}{$\uparrow 1.68 \pm 2.52\%$}}} & $34.05 \pm 11.41$ \\[1pt]
\hline
\end{tabular}

\par\smallskip
\begin{minipage}{\textwidth}
\footnotesize
\textit{Notes.} Entries are means $\pm$ sample SDs over ten instances per family. Average pools all 60 instances, giving each family equal weight. Revenue is in \$1,000. Bold marks the highest average revenue among the eight online methods. Both positive-regularization LSOA columns use $\theta=10$; the calibrated column uses $D^{\mathrm{val}}=2$. Colored lines show paired percentage changes from Greedy, with SDs in percentage points. Green $\uparrow$ denotes a gain and red $\downarrow$ a loss. These conventions also apply to the New York revenue table.
\end{minipage}
\end{table}

LSOA outperforms Greedy in all six demand scenarios. Mean paired gains range from 0.69\% under Peaks to 30.61\% under Commute, with a 10.99\% gain under Hotspot. The larger gains in Hotspot and Commute suggest that reserving vehicles for future requests is particularly valuable when demand is spatially concentrated or its main direction changes during the day. The family-level leaders are Batch in Poisson, Peaks, Overdispersed, and Drift, LSOA ($\theta=0$) in Hotspot, and Sampled Lookahead in Commute.

Across all 60 instances, LSOA has the highest average revenue among the eight online methods at \$30.07k, followed by LSOA ($\theta=0$) at \$29.86k and Sampled Lookahead at \$29.82k. These averages give equal weight to each demand family. Together with the positive gains in all six families, they show that LSOA maintains good performance across the different demand patterns, while the best method within a particular family varies.

Relative to LSOA ($\gamma=1/2$), calibration raises mean paired revenue by 4.93\% across all 60 instances. Gains range from 4.99\% to 7.90\% in five families, while Commute selects $\gamma=0.45$ and gains 0.65\%. The comparison shows why the threshold should be calibrated to the demand pattern: a lower cutoff helps in several families, but is not uniformly preferable. These gains do not require changing the mean-demand planning model or improving its forecasts.

Mean Demand Lookahead earns less than Greedy under Poisson, Peaks, Overdispersed, and Drift. Its admission test compares the current reward $r_{\xi_i}$ with the loss in forecast revenue $W_i^{\mathrm{LP}}(0)-W_i^{\mathrm{LP}}(1)$ from committing a whole vehicle. The continuation LP allocates vehicles to expected demand counts as if those counts were available for service. This can make waiting appear more profitable than it is on the realized request stream. Under Poisson demand, Mean Demand Lookahead serves 61.23\% of requests, compared with 67.06\% for Greedy, which is consistent with reserving too much capacity for forecast demand. LSOA instead uses prices from a plan that includes the current request and applies a validation-selected admission threshold. Mean Demand Lookahead's largest gain over Greedy occurs under Commute, where the benefit of reserving vehicles is larger.

\begin{table}[!ht]
\centering
\caption{Synthetic computation time.}\label{tab:synthetic-time}
\fontsize{8}{10}\selectfont
\setlength{\tabcolsep}{2pt}
\renewcommand{\arraystretch}{1.15}
\begin{tabular}{@{}lcccccccc@{}}
\hline
Demand & Greedy & LSOA & \shortstack{LSOA\\($\gamma=1/2$)} & \shortstack{LSOA\\($\theta=0$)} & \shortstack{Mean Demand\\Lookahead} & \shortstack{Sampled\\Lookahead} & \shortstack{Supply-aware\\MBP} & Batch \\
\hline
Poisson & $0.006 \pm 0.000$ & $6.308 \pm 0.117$ & $8.134 \pm 0.139$ & $1.657 \pm 0.056$ & $3.771 \pm 0.067$ & $36.134 \pm 0.428$ & $0.00524 \pm 0.00006$ & $0.572 \pm 0.010$ \\
Peaks & $0.006 \pm 0.000$ & $5.403 \pm 0.097$ & $7.828 \pm 0.086$ & $1.581 \pm 0.027$ & $3.786 \pm 0.056$ & $38.495 \pm 0.451$ & $0.00508 \pm 0.00003$ & $0.566 \pm 0.003$ \\
Overdispersed & $0.006 \pm 0.000$ & $5.724 \pm 0.459$ & $7.828 \pm 0.487$ & $1.514 \pm 0.089$ & $3.644 \pm 0.166$ & $33.750 \pm 2.364$ & $0.00546 \pm 0.00083$ & $0.576 \pm 0.012$ \\
Hotspot & $0.004 \pm 0.000$ & $4.293 \pm 0.073$ & $6.764 \pm 0.138$ & $1.093 \pm 0.036$ & $3.556 \pm 0.078$ & $35.138 \pm 0.609$ & $0.00440 \pm 0.00090$ & $0.551 \pm 0.006$ \\
Commute & $0.003 \pm 0.000$ & $6.612 \pm 0.378$ & $6.296 \pm 0.242$ & $0.691 \pm 0.027$ & $2.853 \pm 0.097$ & $33.239 \pm 0.743$ & $0.00371 \pm 0.00005$ & $0.545 \pm 0.013$ \\
Drift & $0.006 \pm 0.000$ & $5.520 \pm 0.100$ & $7.804 \pm 0.112$ & $1.475 \pm 0.051$ & $3.608 \pm 0.032$ & $31.727 \pm 0.586$ & $0.00533 \pm 0.00076$ & $0.598 \pm 0.015$ \\
\hline
Average & $0.005 \pm 0.001$ & $5.643 \pm 0.785$ & $7.442 \pm 0.713$ & $1.335 \pm 0.345$ & $3.536 \pm 0.331$ & $34.747 \pm 2.437$ & $0.00487 \pm 0.00084$ & $0.568 \pm 0.021$ \\
\hline
\end{tabular}

\par\smallskip
\begin{minipage}{\textwidth}
\footnotesize
\textit{Notes.} Mean $\pm$ sample SD of instance-level mean times, in ms/request (ms/nonempty dispatch epoch for Batch). Each family has ten instances, and Average pools all 60. Each configuration runs once per instance. All requests, including infeasible ones, enter sequential-method timing. The New York time table uses the same timing conventions. Validation replays are performed offline and are excluded from per-request computation time.
\end{minipage}
\end{table}

Table~\ref{tab:synthetic-time} reports mean decision times for each family. LSOA earns more than Sampled Lookahead in five of the six families and takes less time in all six. Under Commute, LSOA's mean paired revenue difference relative to Sampled Lookahead is $-0.37\%$. Increasing the number of sampled paths from 20 to 50 raises mean paired revenue by 0.07\% while increasing average decision time from 34.75 to 87.14 milliseconds per request. Supply-aware MBP needs only the idle fleet and the occupation price and averages 0.00487 milliseconds per request. Its mean paired revenue falls 7.72\% below LSOA across the 60 instances, with the largest shortfall under Commute.

\subsubsection{Sensitivity to the calibration sample.}
\begin{figure}[!ht]
\centering
\includegraphics[width=\textwidth]{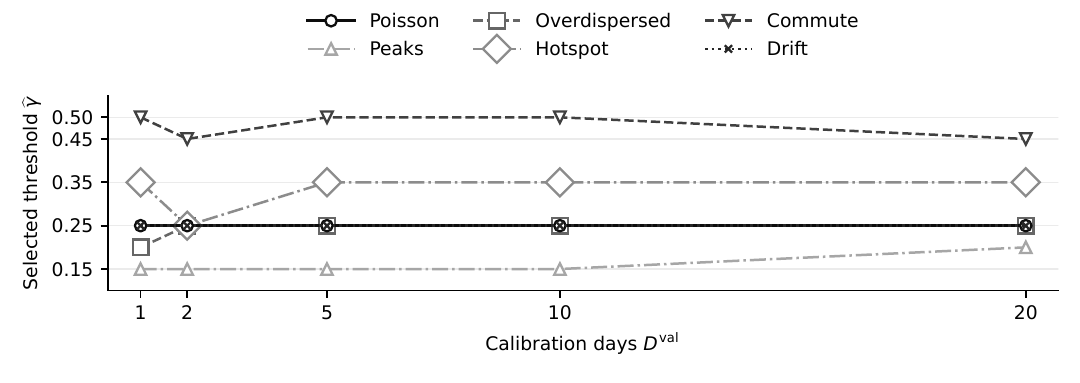}
\caption{Selected thresholds at different calibration sample sizes.}
\label{fig:validation-size}
\end{figure}

We select one threshold per demand family using the first 1, 2, 5, 10, or 20 days of the same validation sample. All five comparisons use the same regularization, forecasting histories, and ten test streams per family. Figure~\ref{fig:validation-size} shows the threshold selected at each sample size.

Poisson and Drift select the same threshold at all five sample sizes. Hotspot selects $\gamma=0.25$ with two days and $\gamma=0.35$ with twenty. Its mean paired test-revenue gain over $\gamma=1/2$ rises from 4.99\% to 6.04\%, the largest increase among the six families. Across all 60 instances, the corresponding gain changes from 4.93\% to 5.08\%. These comparisons measure the effect of adding days to the observed validation sample.
\FloatBarrier

\subsection{New York results}\label{sec:nyc-results}

For New York, we fix $\theta=1$ and $D^{\mathrm{val}}=2$. Validation selects $\gamma=0.05$ for each of the three fleets. We also report fixed $\gamma=1/2$ and LSOA ($\theta=0$), as in the synthetic comparison.

Table~\ref{tab:nyc-revenue} reports daily revenue for the 25-region area. The comparison uses the same ten test dates at each fleet size: February 6--7, 10--14, and 17--19.

\begin{table}[!ht]
\centering
\caption{New York revenue.}\label{tab:nyc-revenue}
\fontsize{7.5}{9}\selectfont
\renewcommand{\scriptsize}{\fontsize{6.8}{8}\selectfont}
\setlength{\tabcolsep}{1.2pt}
\renewcommand{\arraystretch}{1.15}
\begin{tabular}{@{}lccccccccc@{}}
\hline\noalign{\vskip 4pt}
Fleet $C$ & Greedy & LSOA & \shortstack{LSOA\\($\gamma=1/2$)} & \shortstack{LSOA\\($\theta=0$)} & \shortstack{Mean Demand\\Lookahead} & \shortstack{Sampled\\Lookahead} & \shortstack{Supply-aware\\MBP} & Batch & Upper bound \\
\hline\noalign{\vskip 4pt}
250 & $145.41 \pm 5.18$ & \shortstack[c]{$154.96 \pm 7.45$ \\ {\scriptsize \textcolor[rgb]{0.08,0.42,0.23}{$\uparrow 6.53 \pm 1.92\%$}}} & \shortstack[c]{$135.13 \pm 6.71$ \\ {\scriptsize \textcolor[rgb]{0.70,0.15,0.13}{$\downarrow 7.11 \pm 1.82\%$}}} & \shortstack[c]{$149.71 \pm 7.61$ \\ {\scriptsize \textcolor[rgb]{0.08,0.42,0.23}{$\uparrow 2.91 \pm 2.11\%$}}} & \shortstack[c]{$136.37 \pm 7.03$ \\ {\scriptsize \textcolor[rgb]{0.70,0.15,0.13}{$\downarrow 6.26 \pm 2.01\%$}}} & \shortstack[c]{$156.87 \pm 6.98$ \\ {\scriptsize \textcolor[rgb]{0.08,0.42,0.23}{$\uparrow 7.86 \pm 1.58\%$}}} & \shortstack[c]{$147.75 \pm 5.58$ \\ {\scriptsize \textcolor[rgb]{0.08,0.42,0.23}{$\uparrow 1.60 \pm 0.69\%$}}} & \shortstack[c]{$147.57 \pm 4.99$ \\ {\scriptsize \textcolor[rgb]{0.08,0.42,0.23}{$\uparrow 1.50 \pm 1.13\%$}}} & $167.58 \pm 5.83$ \\[1pt]
500 & $220.44 \pm 14.72$ & \shortstack[c]{$265.86 \pm 20.17$ \\ {\scriptsize \textcolor[rgb]{0.08,0.42,0.23}{$\uparrow 20.56 \pm 3.60\%$}}} & \shortstack[c]{$252.64 \pm 18.21$ \\ {\scriptsize \textcolor[rgb]{0.08,0.42,0.23}{$\uparrow 14.59 \pm 3.38\%$}}} & \shortstack[c]{$262.31 \pm 19.73$ \\ {\scriptsize \textcolor[rgb]{0.08,0.42,0.23}{$\uparrow 18.95 \pm 3.56\%$}}} & \shortstack[c]{$252.68 \pm 18.32$ \\ {\scriptsize \textcolor[rgb]{0.08,0.42,0.23}{$\uparrow 14.62 \pm 3.48\%$}}} & \shortstack[c]{$272.17 \pm 21.20$ \\ {\scriptsize \textcolor[rgb]{0.08,0.42,0.23}{$\uparrow 23.40 \pm 3.94\%$}}} & \shortstack[c]{$232.04 \pm 16.66$ \\ {\scriptsize \textcolor[rgb]{0.08,0.42,0.23}{$\uparrow 5.23 \pm 1.26\%$}}} & \shortstack[c]{$218.07 \pm 12.97$ \\ {\scriptsize \textcolor[rgb]{0.70,0.15,0.13}{$\downarrow 1.01 \pm 1.94\%$}}} & $290.84 \pm 22.01$ \\[1pt]
750 & $254.88 \pm 19.57$ & \shortstack[c]{$306.63 \pm 29.01$ \\ {\scriptsize \textcolor[rgb]{0.08,0.42,0.23}{$\uparrow 20.15 \pm 3.94\%$}}} & \shortstack[c]{$302.56 \pm 28.27$ \\ {\scriptsize \textcolor[rgb]{0.08,0.42,0.23}{$\uparrow 18.57 \pm 3.95\%$}}} & \shortstack[c]{$306.36 \pm 28.96$ \\ {\scriptsize \textcolor[rgb]{0.08,0.42,0.23}{$\uparrow 20.04 \pm 3.90\%$}}} & \shortstack[c]{$302.75 \pm 28.16$ \\ {\scriptsize \textcolor[rgb]{0.08,0.42,0.23}{$\uparrow 18.65 \pm 3.88\%$}}} & \shortstack[c]{$314.79 \pm 31.28$ \\ {\scriptsize \textcolor[rgb]{0.08,0.42,0.23}{$\uparrow 23.30 \pm 4.58\%$}}} & \shortstack[c]{$273.53 \pm 23.11$ \\ {\scriptsize \textcolor[rgb]{0.08,0.42,0.23}{$\uparrow 7.25 \pm 1.65\%$}}} & \shortstack[c]{$250.65 \pm 17.99$ \\ {\scriptsize \textcolor[rgb]{0.70,0.15,0.13}{$\downarrow 1.60 \pm 1.67\%$}}} & $333.41 \pm 34.69$ \\[1pt]
\hline\noalign{\vskip 4pt}
Average & $206.91 \pm 48.53$ & \shortstack[c]{$242.48 \pm 68.22$ \\ {\scriptsize \textcolor[rgb]{0.08,0.42,0.23}{$\uparrow 15.74 \pm 7.34\%$}}} & \shortstack[c]{$230.11 \pm 73.90$ \\ {\scriptsize \textcolor[rgb]{0.08,0.42,0.23}{$\uparrow 8.68 \pm 11.88\%$}}} & \shortstack[c]{$239.46 \pm 70.00$ \\ {\scriptsize \textcolor[rgb]{0.08,0.42,0.23}{$\uparrow 13.97 \pm 8.57\%$}}} & \shortstack[c]{$230.60 \pm 73.42$ \\ {\scriptsize \textcolor[rgb]{0.08,0.42,0.23}{$\uparrow 9.00 \pm 11.53\%$}}} & \shortstack[c]{$\boldsymbol{247.95 \pm 71.14}$ \\ {\scriptsize \textcolor[rgb]{0.08,0.42,0.23}{$\boldsymbol{\uparrow 18.19 \pm 8.20\%}$}}} & \shortstack[c]{$217.78 \pm 55.63$ \\ {\scriptsize \textcolor[rgb]{0.08,0.42,0.23}{$\uparrow 4.70 \pm 2.67\%$}}} & \shortstack[c]{$205.43 \pm 45.55$ \\ {\scriptsize \textcolor[rgb]{0.70,0.15,0.13}{$\downarrow 0.37 \pm 2.07\%$}}} & $263.94 \pm 75.17$ \\[1pt]
\hline
\end{tabular}

\par\smallskip
\begin{minipage}{\textwidth}
\footnotesize
\textit{Notes.} Entries are means $\pm$ sample SDs across the same ten test dates for each fleet; revenue is in \$1,000. Average pools the 30 date--fleet settings, giving each fleet size equal weight. Colored lines report paired gains over Greedy, following Table~\ref{tab:synthetic-performance}.
\end{minipage}
\end{table}

LSOA earns more than Greedy and Mean Demand Lookahead at all three fleet sizes. Its mean paired gains over Greedy are 6.53\%, 20.56\%, and 20.15\% at $C=250$, $500$, and $750$, respectively. Sampled Lookahead has the highest mean revenue at each fleet size and averages \$247.95k across the 30 date--fleet settings. LSOA follows at \$242.48k, compared with \$239.46k for LSOA ($\theta=0$).

The calibrated threshold raises mean paired revenue by 7.08\% relative to $\gamma=1/2$. Gains decrease from 14.69\% at $C=250$ to 5.20\% at $C=500$ and 1.33\% at $C=750$, showing that adjusting the mean-demand admission test matters most for the smallest fleet. LSOA also exceeds LSOA ($\theta=0$) by 3.53\%, 1.35\%, and 0.09\% in mean paired revenue at these fleet sizes. It earns more on all ten test dates at $C=250$ and $500$, and on eight dates at $C=750$.

Batch gains 1.50\% over Greedy at $C=250$, but loses 1.01\% and 1.60\% at $C=500$ and $750$. It selects rides by current fare without accounting for subsequent service opportunities, while the average dispatch wait of 1.03 minutes delays vehicle reuse. These delays may offset the benefit of selecting requests jointly.

Supply-aware MBP improves mean paired revenue over Greedy by 1.60--7.25\% across the three fleets. LSOA earns a further 4.84--14.55\% relative to MBP, indicating that current idle supply and an aggregate occupation price alone do not recover the gains from forecast-based planning in these instances.

Sampled Lookahead compares acceptance and rejection by optimizing vehicle use on the same sampled future request streams. The sample networks retain whole vehicles and continuous arrival and completion times, capturing how each action changes subsequent service opportunities. This more detailed continuation model may explain its higher revenue in New York, but it requires two optimizations per sample at each physically feasible arrival.

\begin{table}[!ht]
\centering
\caption{New York computation time.}\label{tab:nyc-time}
\fontsize{8}{10}\selectfont
\setlength{\tabcolsep}{2pt}
\renewcommand{\arraystretch}{1.15}
\begin{tabular}{@{}lcccccccc@{}}
\hline
Fleet $C$ & Greedy & LSOA & \shortstack{LSOA\\($\gamma=1/2$)} & \shortstack{LSOA\\($\theta=0$)} & \shortstack{Mean Demand\\Lookahead} & \shortstack{Sampled\\Lookahead} & \shortstack{Supply-aware\\MBP} & Batch \\
\hline
250 & $0.005 \pm 0.001$ & $26.836 \pm 3.586$ & $37.291 \pm 6.302$ & $4.584 \pm 0.774$ & $10.374 \pm 1.455$ & $308.649 \pm 42.351$ & $0.00508 \pm 0.00024$ & $1.124 \pm 0.077$ \\
500 & $0.009 \pm 0.005$ & $37.874 \pm 7.287$ & $39.683 \pm 6.934$ & $5.530 \pm 1.319$ & $10.917 \pm 1.966$ & $395.945 \pm 77.253$ & $0.00715 \pm 0.00317$ & $1.409 \pm 0.640$ \\
750 & $0.009 \pm 0.003$ & $32.766 \pm 5.518$ & $32.970 \pm 5.102$ & $4.477 \pm 0.913$ & $9.471 \pm 1.736$ & $368.634 \pm 66.248$ & $0.00863 \pm 0.00382$ & $1.526 \pm 0.786$ \\
\hline
Average & $0.008 \pm 0.004$ & $32.492 \pm 7.139$ & $36.648 \pm 6.581$ & $4.864 \pm 1.103$ & $10.254 \pm 1.778$ & $357.743 \pm 71.736$ & $0.00695 \pm 0.00314$ & $1.353 \pm 0.592$ \\
\hline
\end{tabular}

\par\smallskip
\begin{minipage}{\textwidth}
\footnotesize
\textit{Notes.} Means $\pm$ sample SDs across the same ten dates at each fleet size. Average pools the 30 date--fleet settings. Configurations follow Table~\ref{tab:nyc-revenue}. Times are milliseconds per observed request, including infeasible arrivals, or per nonempty dispatch epoch for Batch. Offline validation is excluded.
\end{minipage}
\end{table}

Table~\ref{tab:nyc-time} shows the computation required for these revenues. LSOA averages 32.49 milliseconds per request, compared with 357.74 for Sampled Lookahead. Its mean paired percentage revenue loss relative to Sampled Lookahead ranges from 1.24\% to 2.54\% across the three fleet sizes. LSOA ($\theta=0$) takes 4.86 milliseconds on average, but the revenue difference from LSOA is larger at the smallest fleet. It also earns more than Mean Demand Lookahead at every fleet size and requires less computation time. Supply-aware MBP averages 0.00695 milliseconds per request, while Batch averages 1.35 milliseconds per nonempty dispatch epoch.

\section{Conclusion and Discussion}\label{sec:conclusion}\label{sec:discussion}

We study how an operator with a fixed fleet should accept or reject ride requests when future demand is unknown. Our method, LSOA, uses a rolling mean-demand plan that accounts for where and when vehicles become available after each trip. The plan uses dual prices to estimate the opportunity cost of dispatching a vehicle. The operator accepts a feasible request when its planned service fraction reaches a common threshold, which we calibrate using revenue on validation days.

For one vehicle on a single route, we derive the exact range of thresholds that attain optimal expected revenue. When waiting is optimal, the same threshold can reject the current request and accept the first higher-reward request that arrives. When several origins share return demand at a hub, the fleet may need only some of the available vehicles. If the plan assigns the same fraction to each initial outbound request, a common threshold sends either all vehicles or none. The resulting revenue loss can grow in proportion to fleet size even under the best common threshold. Independent validation bounds the revenue loss from threshold selection relative to the best candidate. At a fixed threshold, under Poisson arrivals and a margin condition, the revenue gap between using estimated and true mean demand has a tight $O(D^{-1/2})$ bound with $D$ historical days.

The New York results distinguish where it pays to calibrate the threshold from where it pays to anticipate later requests. With 250 vehicles, calibration raises revenue by 14.69\% over the zero-cutoff rule, compared with 1.33\% with 750 vehicles. Yet LSOA's gain over serving every feasible ride is 6.53\% with 250 vehicles and about 20\% with either larger fleet. In the synthetic study, gains over Greedy are larger when demand is concentrated or changes direction during the day. A fleet deciding where to use anticipatory admission should therefore examine which locations run short of vehicles and when vehicles return. To calibrate the policy for that market, validation days should reflect its fleet size, directional flows, and time-of-day demand.

For implementation, the fleet must record where and when each busy vehicle will next become available. It must also log incoming requests that it turns away or cannot serve. If only completed trips enter the forecasting history, rejecting an order removes an observation of demand at that origin. The New York replay uses completed taxi trips as the common request stream and holds supply, fares, and trip durations fixed. Estimating how a change in fares affects demand would require records of offers, quoted prices, and passenger responses.

The revenue objective values which trips the fleet can serve over the day. Passengers may also respond to waits and rejections. For example, \citet{liu2026waiting} find that passengers' responses to updated delay announcements depend on how long they have already waited. A fleet pilot could measure service rates, waiting times, and repeat use alongside revenue. Pickup travel and empty relocation are natural extensions because both change which vehicle can reach the next passenger. Extending the calibration guarantee to markets where demand changes between validation and operation remains an open question.

\ACKNOWLEDGMENT{This work was supported by the Tsinghua University Talent Introduction Fund [grant number 53330900126].}

\begingroup
\renewcommand{\baselinestretch}{1.3}\selectfont
\def\bibfont{\fontsize{9}{11}\selectfont}
\def\bibsep{2pt}
\bibliography{ride-hailing-ts}

@article{talluri1998bidprice,
  author = {Talluri, Kalyan and van Ryzin, Garrett},
  title = {An Analysis of Bid-Price Controls for Network Revenue Management},
  journal = {Management Science},
  year = {1998},
  volume = {44},
  number = {11},
  pages = {1577--1593},
  doi = {10.1287/mnsc.44.11.1577}
}

@article{talluri1999randomized,
  author = {Talluri, Kalyan and van Ryzin, Garrett},
  title = {A Randomized Linear Programming Method for Computing Network Bid Prices},
  journal = {Transportation Science},
  year = {1999},
  volume = {33},
  number = {2},
  pages = {207--216},
  doi = {10.1287/trsc.33.2.207}
}

@article{bimpikis2019spatial,
  author = {Bimpikis, Kostas and Candogan, Ozan and Saban, Daniela},
  title = {Spatial Pricing in Ride-Sharing Networks},
  journal = {Operations Research},
  year = {2019},
  volume = {67},
  number = {3},
  pages = {744--769},
  doi = {10.1287/opre.2018.1800}
}

@article{braverman2019empty,
  author = {Braverman, Anton and Dai, J. G. and Liu, Xin and Ying, Lei},
  title = {Empty-Car Routing in Ridesharing Systems},
  journal = {Operations Research},
  year = {2019},
  volume = {67},
  number = {5},
  pages = {1437--1452},
  doi = {10.1287/opre.2018.1822}
}

@article{ozkan2020matching,
  author = {{\"O}zkan, Erhun and Ward, Amy R.},
  title = {Dynamic Matching for Real-Time Ride Sharing},
  journal = {Stochastic Systems},
  year = {2020},
  volume = {10},
  number = {1},
  pages = {29--70},
  doi = {10.1287/stsy.2019.0037}
}

@article{feng2021way,
  author = {Feng, Guiyun and Kong, Guangwen and Wang, Zizhuo},
  title = {We Are on the Way: Analysis of On-Demand Ride-Hailing Systems},
  journal = {Manufacturing \& Service Operations Management},
  year = {2021},
  volume = {23},
  number = {5},
  pages = {1237--1256},
  doi = {10.1287/msom.2020.0880}
}

@article{gong2022assortment,
  author = {Gong, Xiao-Yue and Goyal, Vineet and Iyengar, Garud N. and Simchi-Levi, David and Udwani, Rajan and Wang, Shuangyu},
  title = {Online Assortment Optimization with Reusable Resources},
  journal = {Management Science},
  year = {2022},
  volume = {68},
  number = {7},
  pages = {4772--4785},
  doi = {10.1287/mnsc.2021.4134}
}

@article{goyal2025reusable,
  author = {Goyal, Vineet and Iyengar, Garud and Udwani, Rajan},
  title = {Asymptotically Optimal Competitive Ratio for Online Allocation of Reusable Resources},
  journal = {Operations Research},
  year = {2025},
  volume = {73},
  number = {4},
  pages = {1897--1915},
  doi = {10.1287/opre.2021.0695}
}

@article{balseiro2021relocating,
  author = {Balseiro, Santiago R. and Brown, David B. and Chen, Chen},
  title = {Dynamic Pricing of Relocating Resources in Large Networks},
  journal = {Management Science},
  year = {2021},
  volume = {67},
  number = {7},
  pages = {4075--4094},
  doi = {10.1287/mnsc.2020.3735}
}

@article{balseiro2023mirror,
  author = {Balseiro, Santiago R. and Lu, Haihao and Mirrokni, Vahab},
  title = {The Best of Many Worlds: Dual Mirror Descent for Online Allocation Problems},
  journal = {Operations Research},
  year = {2023},
  volume = {71},
  number = {1},
  pages = {101--119},
  doi = {10.1287/opre.2021.2242}
}

@article{jasin2012resolving,
  author = {Jasin, Stefanus and Kumar, Sunil},
  title = {A Re-Solving Heuristic with Bounded Revenue Loss for Network Revenue Management with Customer Choice},
  journal = {Mathematics of Operations Research},
  year = {2012},
  volume = {37},
  number = {2},
  pages = {313--345},
  doi = {10.1287/moor.1120.0537}
}

@article{brown2010information,
  author = {Brown, David B. and Smith, James E. and Sun, Peng},
  title = {Information Relaxations and Duality in Stochastic Dynamic Programs},
  journal = {Operations Research},
  year = {2010},
  volume = {58},
  number = {4-part-1},
  pages = {785--801},
  doi = {10.1287/opre.1090.0796}
}

@article{bertsekas1999rollout,
  author = {Bertsekas, Dimitri P. and Castanon, David A.},
  title = {Rollout Algorithms for Stochastic Scheduling Problems},
  journal = {Journal of Heuristics},
  year = {1999},
  volume = {5},
  number = {1},
  pages = {89--108},
  doi = {10.1023/A:1009634810396}
}

@article{bumpensanti2020resolving,
  author = {Bumpensanti, Pornpawee and Wang, He},
  title = {A Re-Solving Heuristic with Uniformly Bounded Loss for Network Revenue Management},
  journal = {Management Science},
  year = {2020},
  volume = {66},
  number = {7},
  pages = {2993--3009},
  doi = {10.1287/mnsc.2019.3365}
}

@misc{nyctlc2025data,
  author = {{New York City Taxi and Limousine Commission}},
  title = {{TLC} trip record data: Yellow taxi trips, {January} and {February} 2025},
  year = {2025},
  howpublished = {\url{https://www.nyc.gov/site/tlc/about/tlc-trip-record-data.page}},
  note = {Taxi-zone lookup and yellow-taxi data dictionary. Accessed September 16, 2026}
}

@article{akan2009bidprice,
  author = {Akan, Mustafa and Ata, Bar{\i}{\c{s}}},
  title = {Bid-Price Controls for Network Revenue Management: Martingale Characterization of Optimal Bid Prices},
  journal = {Mathematics of Operations Research},
  year = {2009},
  volume = {34},
  number = {4},
  pages = {912--936},
  doi = {10.1287/moor.1090.0411}
}

@article{topaloglu2008stochastic,
  author = {Topaloglu, Huseyin},
  title = {A Stochastic Approximation Method to Compute Bid Prices in Network Revenue Management Problems},
  journal = {INFORMS Journal on Computing},
  year = {2008},
  volume = {20},
  number = {4},
  pages = {596--610},
  doi = {10.1287/ijoc.1080.0269}
}

@article{topaloglu2007sensitivity,
  author = {Topaloglu, Huseyin and Powell, Warren B.},
  title = {Sensitivity Analysis of a Dynamic Fleet Management Model Using Approximate Dynamic Programming},
  journal = {Operations Research},
  year = {2007},
  volume = {55},
  number = {2},
  pages = {319--331},
  doi = {10.1287/opre.1060.0347}
}

@article{li2023bidprices,
  author = {Li, Dong and Pang, Zhan and Qian, Lixian},
  title = {Bid price controls for car rental network revenue management},
  journal = {Production and Operations Management},
  year = {2023},
  volume = {32},
  number = {1},
  pages = {261--282},
  doi = {10.1111/poms.13836}
}

@article{chen2024spatial,
  author = {Chen, Qi and Lei, Yanzhe and Jasin, Stefanus},
  title = {Real-Time Spatial--Intertemporal Pricing and Relocation in a Ride-Hailing Network: Near-Optimal Policies and the Value of Dynamic Pricing},
  journal = {Operations Research},
  year = {2024},
  volume = {72},
  number = {5},
  pages = {2097--2118},
  doi = {10.1287/opre.2022.2425}
}

@inproceedings{xu2018dispatch,
  author = {Xu, Zhe and Li, Zhixin and Guan, Qingwen and Zhang, Dingshui and Li, Qiang and Nan, Junxiao and Liu, Chunyang and Bian, Wei and Ye, Jieping},
  title = {Large-Scale Order Dispatch in On-Demand Ride-Hailing Platforms: A Learning and Planning Approach},
  booktitle = {Proceedings of the 24th ACM SIGKDD International Conference on Knowledge Discovery and Data Mining},
  year = {2018},
  pages = {905--913},
  doi = {10.1145/3219819.3219824}
}

@inproceedings{dickerson2018allocation,
  author = {Dickerson, John P. and Sankararaman, Karthik A. and Srinivasan, Aravind and Xu, Pan},
  title = {Allocation Problems in Ride-Sharing Platforms: Online Matching with Offline Reusable Resources},
  booktitle = {Proceedings of the AAAI Conference on Artificial Intelligence},
  year = {2018},
  volume = {32},
  pages = {1007--1014},
  doi = {10.1609/aaai.v32i1.11477},
  url = {https://jpdickerson.com/pubs/dickerson18allocation.pdf}
}

@article{brown2025fluid,
  author = {Brown, David B. and Zhang, Jingwei},
  title = {Fluid Policies, Reoptimization, and Performance Guarantees in Dynamic Resource Allocation},
  journal = {Operations Research},
  year = {2025},
  volume = {73},
  number = {2},
  pages = {1029--1045},
  doi = {10.1287/opre.2022.0601}
}

@article{vera2020prophet,
  author = {Vera, Alberto and Banerjee, Siddhartha},
  title = {The Bayesian Prophet: A Low-Regret Framework for Online Decision Making},
  journal = {Management Science},
  year = {2021},
  volume = {67},
  number = {3},
  pages = {1368--1391},
  doi = {10.1287/mnsc.2020.3624}
}

@article{kanoria2022mirror,
  author = {Kanoria, Yash and Qian, Pengyu},
  title = {Blind Dynamic Resource Allocation in Closed Networks via Mirror Backpressure},
  year = {2024},
  journal = {Management Science},
  volume = {70},
  number = {8},
  pages = {5445--5462},
  doi = {10.1287/mnsc.2023.4934}
}

@article{zhang2025reusable,
  author = {Zhang, Xilin and Cheung, Wang Chi},
  title = {Online Allocation of Reusable Resources in Nonstationary Environments},
  journal = {Mathematics of Operations Research},
  year = {2026},
  volume = {51},
  number = {3},
  pages = {1853--1883},
  doi = {10.1287/moor.2023.0250}
}

@article{afeche2023networks,
  author = {Af\`eche, Philipp and Liu, Zhe and Maglaras, Costis},
  title = {Ride-Hailing Networks with Strategic Drivers: The Impact of Platform Control Capabilities on Performance},
  journal = {Manufacturing \& Service Operations Management},
  year = {2023},
  volume = {25},
  number = {5},
  pages = {1890--1908},
  doi = {10.1287/msom.2023.1221}
}

@book{boyd2004convex,
  author = {Boyd, Stephen and Vandenberghe, Lieven},
  title = {Convex Optimization},
  publisher = {Cambridge University Press},
  year = {2004},
  url = {https://web.stanford.edu/~boyd/cvxbook/bv_cvxbook.pdf}
}

@article{liu2026waiting,
  author = {Liu, Yuhan and Qi, Jin and Yang, Hai and Tang, Yili},
  title = {Passenger Waiting Behavior with Sunk Waiting Time: An Empirical Study in the Ride-Sourcing Market},
  journal = {Production and Operations Management},
  year = {2026},
  note = {OnlineFirst},
  doi = {10.1177/10591478261488817}
}

@misc{blacklane2026partners,
  author = {{Blacklane}},
  title = {What does it mean to become a {Blacklane} partner?},
  year = {2026},
  howpublished = {Blacklane Partner Help Center, \url{https://partner-help.blacklane.com/en/articles/8420677-what-does-it-mean-to-become-a-blacklane-partner}},
  note = {Accessed November 21, 2026}
}

@misc{bookingTaxiSupplier,
  author = {{Booking{\char46}com}},
  title = {{Booking.com Taxi Supplier API} documentation},
  year = {2026},
  howpublished = {\url{https://connect.taxi.booking.com/}},
  note = {Accessed November 21, 2026}
}

@inproceedings{kakade2002approximate,
  author = {Kakade, Sham and Langford, John},
  title = {Approximately Optimal Approximate Reinforcement Learning},
  booktitle = {Proceedings of the Nineteenth International Conference on Machine Learning},
  year = {2002},
  pages = {267--274},
  publisher = {Morgan Kaufmann}
}
\endgroup
\ECRUNAUTHOR{}
\ECAUpunct{}
\ECSwitch
% Keep appendix hyperlink destinations distinct from the main paper.
\renewcommand{\theHsection}{EC.\arabic{section}}
\renewcommand{\theHtable}{EC.\arabic{table}}
\renewcommand{\theHequation}{EC.\arabic{equation}}
\renewcommand{\theHlemma}{EC.\arabic{lemma}}
\renewcommand{\theHtheorem}{EC.\arabic{theorem}}
\renewcommand{\theHcorollary}{EC.\arabic{corollary}}
\renewcommand{\theHproposition}{EC.\arabic{proposition}}
\renewcommand{\theHassumption}{EC.\arabic{assumption}}
\ECHead{Electronic Companion}
\ECHead{Learned spatial opportunity cost admission for ride-hailing with calibrated thresholds}
\section{Vehicle plans and solver prices}

\subsection{Proof of Proposition~\ref{prop:ts-duality}}\label{ec:planning}

\begin{proof}
Fix request $i$, the fleet state, and the forecast counts. Construct the network on $\mathcal N_i=(\mathcal S\times\mathcal T_i)\cup\{v_{\mathrm{term}}\}$. Node $(s,\tau)$ receives supply $b_{s,\tau|i}-b_{s,\tau-1|i}$, with $b_{s,-1|i}=0$. A service arc runs from departure to completion, or to the terminal node if completion is beyond the last grid point. The current arc has capacity one and each forecast arc has capacity $\widehat q_{z,\tau|i}$. Waiting arcs carry unused capacity forward one interval at the same location, and from the last grid point to the terminal node. For a feasible solution of \eqref{eq:quadratic-plan}, let $w_{s,\tau|i}\ge0$ be the slack in \eqref{eq:quadratic-vehicle}, and set $w_{s,-1|i}:=0$. Taking differences between consecutive cumulative constraints gives
\begin{equation}\label{eq:sparse-vehicle-balance}
\begin{aligned}
&w_{s,\tau|i}-w_{s,\tau-1|i} +\mathbbm{1}\{s=o_{\xi_i},\tau=0\}y_i -\mathbbm{1}\{s=d_{\xi_i},\tau=\ell_{\xi_i}\}y_i\\
&\quad+\sum_{z\in\mathcal Z_{\tau|i}:o_z=s}x_{z,\tau|i} -\sum_{\tau'=0}^{\tau-1} \sum_{\substack{z\in\mathcal Z_{\tau'|i}\\d_z=s,\ \tau'+\ell_z=\tau}} x_{z,\tau'|i} =b_{s,\tau|i}-b_{s,\tau-1|i}, \quad \forall s\in\mathcal S,\ \tau\in\mathcal T_i .
\end{aligned}
\end{equation}
Here $\tau'$ indexes earlier departure intervals; the last sum counts rides completing in region $s$ at interval $\tau$. Equation~\eqref{eq:sparse-vehicle-balance} is flow conservation at $(s,\tau)$. Conversely, summing \eqref{eq:sparse-vehicle-balance} over $\tau=0,\ldots,\tau'$ and using $w_{s,\tau'|i}\ge0$ recovers \eqref{eq:quadratic-vehicle}. Thus the feasible plans are exactly the service components of feasible flows in this network. Since $\ell_z\ge1$, every nonterminal arc advances time. The network is acyclic, so each feasible flow decomposes into vehicle paths that respect location and grid completion times. The service constraints form a compact convex polytope containing the zero plan. The continuous objective attains a maximum there. Its Hessian in the service variables is diagonal with entries $-2\theta$ for the current request and $-2\theta/\widehat q_{z,\tau|i}$ for each positive forecast count. Strict concavity gives a unique service plan.
\end{proof}

\subsection{Sparse implementation of the primal plan}\label{ec:sparse-primal}

We minimize $-F_i$ using the flow balances in~\eqref{eq:sparse-vehicle-balance}. Let $\mu_{s,\tau|i}$ multiply the left side minus the right side of the balance at $(s,\tau)$, with $\mu_{s,T_i|i}=0$. For the waiting flow, stationarity gives the multiplier of its nonnegativity constraint as $\mu_{s,\tau|i}-\mu_{s,\tau+1|i}\ge0$. Since this flow is the slack in \eqref{eq:quadratic-vehicle}, the difference is the corresponding cumulative constraint multiplier $\lambda_{s,\tau|i}$. Summing over the remaining intervals in~\eqref{eq:vehicle-price} gives $\rho_{s,\tau|i}=\mu_{s,\tau|i}$.
\subsection{Combining parallel service arcs}\label{ec:parallel-arcs}

The current arc and its same-type forecast arc in interval zero have the same endpoints and fare. For a fixed total $x_i^{\mathrm{tot}}=y_i+x_{\xi_i,0|i}$ and positive count $q=\widehat q_{\xi_i,0|i}$, their combined penalty depends on $y_i^2+x_{\xi_i,0|i}^2/q$. This quantity equals $(x_i^{\mathrm{tot}})^2/(1+q)+(qy_i-x_{\xi_i,0|i})^2/[q(1+q)]$. Thus one arc of capacity $1+q$ gives the same optimal objective, with $y_i^*=x_i^{\mathrm{tot},*}/(1+q)$ for $\theta>0$. This split satisfies both original bounds. At $\theta=0$ it is optimal but need not be unique. At $q=0$ the current arc stands alone.

\subsection{Proof of Proposition~\ref{prop:ts-admission}}\label{ec:price-rule}

\begin{proof}
Fix optimal vehicle multipliers and keep the bounds $0\le y_i\le1$ explicit. By the KKT conditions \citep{boyd2004convex}, $y_i^*$ maximizes $g_i y_i+\theta y_i(1-y_i)$ on $[0,1]$. Since $\theta>0$, projecting the unconstrained maximizer onto this interval gives $y_i^*=\min\{1,\max\{0,(g_i+\theta)/(2\theta)\}\}$. For $0<\gamma<1$, $y_i^*\ge\gamma$ therefore holds exactly when $g_i\ge\theta(2\gamma-1)$. Uniqueness of $y_i^*$ implies that all optimal vehicle multipliers give the same decision.
\end{proof}

\section{Threshold calibration}\label{ec:calibration}

\subsection{Proof of Proposition~\ref{prop:calibration-value}}
\label{ec:calibration-value}

\begin{proof}
Serving the initial request earns $r_{\mathrm{low}}$. Waiting earns $r_{\mathrm{high}}$ if at least one later request arrives, so its optimal expected revenue is $w=r_{\mathrm{high}}(1-e^{-q})$. In the initial plan, let $x$ be the amount reserved for later rides. Both rewards exceed $\theta$, so the plan uses the whole vehicle and assigns $1-x$ to the initial ride. For $0\le x\le q$, the derivative of the objective with respect to $x$ is $r_{\mathrm{high}}-r_{\mathrm{low}}+2\theta-2\theta(1+1/q)x$. Its minimum is $r_{\mathrm{high}}-r_{\mathrm{low}}-2\theta q>0$. Thus $x=q$ and the fleet accepts the initial ride exactly when $\gamma\le1-q$.

If the fleet waits, the next observed ride and the later rides in the forecast have the same reward and use the same vehicle. With remaining mean demand $m$, Section~\ref{ec:parallel-arcs} gives the observed ride a fraction $1/(1+m)$. Since $m\le q$, every $\gamma\le1/(1+q)$ accepts the first later ride. If $\gamma>1/(1+q)$, there is an early interval of positive demand during which $m>\gamma^{-1}-1$. With positive probability the only later ride arrives in that interval and is rejected. Thus the fleet earns $w$ after waiting exactly when $\gamma\le1/(1+q)$. Since $1-q<1/(1+q)$, comparison of $r_{\mathrm{low}}$ and $w$ gives all three intervals in the proposition.
\end{proof}

\noindent\textbf{Extension to $C$ vehicles.} On the same route, let $C\ge1$, $q>0$, $N\sim\operatorname{Poisson}(q)$, and $r_{\mathrm{high}}>r_{\mathrm{low}}>\theta$. If the first action is $a\in\{0,1\}$, accepting later rides until capacity runs out earns $a r_{\mathrm{low}}+r_{\mathrm{high}}\mathbb E\min\{N,C-a\}$. Thus the initial decision compares $r_{\mathrm{low}}$ with $r_{\mathrm{high}}\mathbb P(N\ge C)$, the value of saving the last vehicle. The initial plan maximizes $r_{\mathrm{low}}y+\theta y(1-y)+r_{\mathrm{high}}x+\theta x(1-x/q)$ subject to $0\le y\le1$, $0\le x\le q$, and $x+y\le C$. If capacity binds, substituting $x=C-y$ and imposing the bounds on $y$ gives the expression below. If $C\ge1+q$, both service amounts reach their upper bounds and the same expression gives one.
\[
y_1^{\mathbb P,*}=\min\left\{1,\max\left\{0,C-q,\frac{C-q(r_{\mathrm{high}}-r_{\mathrm{low}})/(2\theta)}{1+q}\right\}\right\}.
\]
A threshold is optimal exactly when its first action $a_1=\mathbbm{1}\{y_1^{\mathbb P,*}\ge\gamma\}$ attains the better of the two initial values and, if $C-a_1>0$, $\gamma\le1/(1+q)$. Indeed, with $b\ge1$ vehicles left the observed fraction is $\min\{1,b/(1+m)\}$, so this upper bound suffices. For necessity, if $\gamma$ exceeds it, there is positive probability that all $C-a_1$ remaining vehicles could serve one ride each in the early interval described above and no rides arrive afterward. Serving them all would require accepting the last one with a single vehicle, which the threshold rejects.

\subsection{A bound for general networks}\label{ec:network-threshold-loss}
Let $\mathcal F_\gamma$ be the feasible arrivals along $\pi_{\theta,\gamma}^{\mathbb P}$. At such an arrival, the expected gain from accepting, given optimal dispatch afterward, is $m_i:=r_{\xi_i}+V_i^*(\boldsymbol s_i^{\mathrm{accept}})-V_i^*(\boldsymbol s_i^{\mathrm{reject}})$. Write $(x)_+=\max\{x,0\}$. Rejecting forgoes $(m_i)_+$ and accepting forgoes $(-m_i)_+$. Summing these losses along the threshold policy gives the performance difference identity \citep{kakade2002approximate}.
\begin{equation}\label{eq:threshold-network-loss}
V^*(\mathbb P)-V_{\pi_{\theta,\gamma}^{\mathbb P}}(\mathbb P)
=\mathbb E\sum_{i\in\mathcal F_\gamma}\left[(m_i)_+\mathbbm{1}\{y_i^{\mathbb P,*}<\gamma\}+(-m_i)_+\mathbbm{1}\{y_i^{\mathbb P,*}\ge\gamma\}\right].
\end{equation}
At each arrival, the loss is the optimal value before the decision minus the current reward and the conditional expected optimal value at the next arrival. Infeasible arrivals contribute zero. Taking expectations cancels successive continuation values, and the terminal value is zero. For an unbounded number of arrivals the same argument follows by truncation, since the expected reward after request $n$ is at most $R\mathbb E[(I-n)_+]\to0$. Each summand is nonnegative, so the loss is zero exactly when the policy takes an optimal action almost surely at every feasible arrival it encounters.

If $|g_i-\theta(2\gamma-1)-m_i|\le\varepsilon$, a strictly suboptimal decision requires $|m_i|\le\varepsilon$, since the cutoff-adjusted reward and the true gain must lie on opposite sides of zero or one must equal zero.

Fix a feasible nonanticipating reference policy with actions $a_i^\circ$. On its path, compute $y_i^*$ using the same forecasts and regularization as the threshold policies. Starting from $\Gamma_0=(0,1)$, intersect $\Gamma_{i-1}$ with $(0,y_i^*]$ after an acceptance and with $(y_i^*,1)$ after a feasible rejection. An infeasible arrival leaves the interval unchanged. Thus $\Gamma_i$ contains the thresholds that agree with the reference through request $i$.

Given the history at a feasible request $i$, let $L_i^a$ and $U_i^a$ bound the expected revenue from action $a$ and optimal continuation, including the current reward. A feasible continuation gives a lower bound and perfect foresight gives an upper bound, both in expectation. Set $d_i=(L_i^{a_i^\circ}-U_i^{1-a_i^\circ})_+$ at feasible requests and $d_i=0$ otherwise. Then
\begin{equation}\label{eq:threshold-prefix-loss}
V^*(\mathbb P)-\sup_{0<\gamma<1}V_{\pi_{\theta,\gamma}^{\mathbb P}}(\mathbb P)
\ \ge\ \inf_{0<\gamma<1}\mathbb E\sum_{i=1}^I d_i\,\mathbbm{1}\{\gamma\in\Gamma_{i-1}\setminus\Gamma_i\}.
\end{equation}
In particular, suppose that with probability at least $p$ the interval becomes empty and every step that removes thresholds up to that point has $d_i\ge d>0$. Then every common threshold loses at least $pd$ in expected revenue.

\begin{proof}
Use the same request stream and fleet assignment convention for both policies, with independent randomization for a randomized reference. Until their first disagreement, the fleet states and planned fractions coincide. The first disagreement therefore occurs at $i$ exactly when $\gamma\in\Gamma_{i-1}\setminus\Gamma_i$. At that state the reference action is worth at least $L_i^{a_i^\circ}$ and the other action at most $U_i^{1-a_i^\circ}$ under optimal continuation. The corresponding summand in \eqref{eq:threshold-network-loss} is at least $d_i$. All other summands are nonnegative. Taking expectations and the infimum over $\gamma$ gives \eqref{eq:threshold-prefix-loss}. On the empty-interval event each threshold incurs at least $d$ at its first disagreement, giving the bound $pd$.
\end{proof}

\subsection{Proof of Proposition~\ref{prop:hub-threshold}}\label{ec:hub-threshold}

If $n$ vehicles are sent initially, the highest expected revenue is $C r_{\mathrm{high}}(1-e^{-q})-n\delta+r_{\mathrm{ret}}\min\{n,k\}$. This expression increases up to $n=k$ and decreases thereafter. As shown below for arbitrary origins, every initial fraction equals $1-q$ in the identical-origin case. A threshold of $1/2$ sends all vehicles and serves all $k$ return rides. A threshold in $(1-q,1/(1+q)]$ keeps all vehicles at their origins and accepts the first later ride at each origin where one arrives. The best threshold therefore attains the better of $n=0$ and $n=C$. Subtracting each revenue from that at $n=k$ gives Proposition~\ref{prop:hub-threshold}.

\noindent\textbf{Different origins and random return demand.} Consider $C$ origins with one idle vehicle at each and an initially empty hub. At time zero, one request from each origin $s$ to the hub pays $r_{s,0}$. Additional requests on the same route pay $r_{s,1}$ and arrive independently across origins as Poisson processes with mean counts $q_s\in(0,1/2)$ on $(0,\Delta)$ and constant rates on this interval. At time $\Delta$, the hub receives $N^{\mathrm{ret}}$ return requests, each paying $r_{\mathrm{ret}}$ for a trip to an origin. The count has any distribution on $\{0,\ldots,C\}$ independent of the outbound requests. All trips last $\Delta$, the horizon is $2\Delta$, and there are no other requests. Simultaneous requests follow the priority order in Section~\ref{sec:problem}. The plan covers the remaining horizon and uses the true conditional mean counts.

A vehicle sent at time zero reaches the hub in time for a return trip. A vehicle held at origin $s$ can earn expected revenue $v_s=r_{s,1}(1-e^{-q_s})$ by serving the first subsequent request there. Let $\delta_s=v_s-r_{s,0}$ be the expected outbound revenue forgone by sending this vehicle now. For a set $\mathcal A$ of origins served at time zero, the maximum expected revenue with those initial decisions is
\begin{equation}\label{eq:hub-revenue}
V(\mathcal A)=\sum_{s=1}^C v_s-\sum_{s\in\mathcal A}\delta_s +r_{\mathrm{ret}}\,\mathbb E\min\{|\mathcal A|,N^{\mathrm{ret}}\}.
\end{equation}
Assume $r_{s,0},r_{\mathrm{ret}}>\theta$ and $r_{s,1}-r_{s,0}>2\theta q_s$ at each origin.

Under the conditions above, let $\mathcal A_k^*$ contain the $k$ origins with the smallest $\delta_s$, breaking ties arbitrarily, and let $\mathcal A_\gamma=\{s:1-q_s\ge\gamma\}$. The optimal online revenue is $V^*(\mathbb P)=\max_{0\le k\le C}V(\mathcal A_k^*)$, and
\begin{equation}\label{eq:hub-gap}
V^*(\mathbb P)-\sup_{0<\gamma<1}V_{\pi_{\theta,\gamma}^{\mathbb P}}(\mathbb P) \ \ge\ \max_{0\le k\le C}V(\mathcal A_k^*)- \max_{0<\gamma<1}V(\mathcal A_\gamma).
\end{equation}
Equality holds when all origins have the same $q_s=q$. We prove both claims below.

For the $k$th vehicle, the expected return reward $r_{\mathrm{ret}}\mathbb P(N^{\mathrm{ret}}\ge k)$ decreases weakly with $k$, while the $k$th smallest outbound cost $\delta_{(k)}$ increases weakly. The operator can dispatch in order of cost until the next cost exceeds its return reward. A common threshold instead orders origins by $1-q_s$. The plan puts initial and later outbound rides in the same departure interval, although only the initial rides reach the hub in time.
\begin{proof}
Fix the set $\mathcal A$ of initial rides served. Every remaining vehicle can serve at most one subsequent outbound ride. Taking the first such ride at each unserved origin earns expected revenue $v_s$ there. Only the $|\mathcal A|$ vehicles dispatched initially are at the hub at $\Delta$. Serving return requests in arrival order until these vehicles or the requests are exhausted earns the last term of \eqref{eq:hub-revenue}. Both choices are feasible online, and each attains its corresponding upper bound. The initial requests are given and reveal no information about later demand. Randomizing their acceptance cannot improve upon the best fixed set. At fixed cardinality $k$, minimizing the second sum in \eqref{eq:hub-revenue} gives $\mathcal A_k^*$.

Next consider the plan at any initial request from origin $s$. Earlier decisions at other origins leave its one vehicle idle. Write $y$ for the current fraction and $x$ for forecast outbound service from $s$. The local constraints are $x\le q_s$ and $x+y\le1$, with $x,y\ge0$. The plan places both flows at time zero, so both reach the hub at $\Delta$. Increasing $y$ with unused local capacity is feasible because additional hub capacity can be left idle, and its marginal objective is at least $r_{s,0}-\theta>0$. Thus $y+x=1$. At fixed total flow the split has no effect on feasibility or revenue elsewhere. As in the proof of Proposition~\ref{prop:calibration-value}, replacing current service by forecast service increases the objective at a rate of at least $r_{s,1}-r_{s,0}-2\theta q_s>0$. Hence $x=q_s$ and $y=1-q_s$ regardless of earlier initial decisions. The threshold sends precisely the vehicles in $\mathcal A_\gamma$, and its continuation revenue is bounded by the corresponding terms of \eqref{eq:hub-revenue}. This proves \eqref{eq:hub-gap}.

If $q_s=q$ throughout, there are only two possible initial sets. With $\gamma=1/2$ the policy sends all $C$ vehicles. At each return request the remaining vehicles suffice for the current request and all remaining forecast demand, since $N^{\mathrm{ret}}\le C$ on every stream. All these rides finish by $H$. Their marginal objectives are positive throughout their feasible ranges because $r_{\mathrm{ret}}>\theta$. The plan fully serves each observed return request, attaining $V(\{1,\ldots,C\})$. A threshold in $(1-q,1/(1+q)]$ sends no vehicle initially. At an outbound arrival $t_i\in(0,\Delta)$, the rolling plan places the return requests at the start of its current interval, when the hub has no idle vehicle. All remaining outbound demand at the current origin lies in this same interval. Its planned fraction is therefore $1/[1+q(1-t_i/\Delta)]$, as in the single-route evaluation in Section~\ref{ec:calibration-value}. The first high-reward request at each origin is served, attaining $V(\varnothing)$. Thus the better of the two initial sets attains the bound.

\end{proof}

\section{Proof of Lemma~\ref{lem:calibration-forecast-support}}
\label{ec:sensitivity-proof}

\begin{proof}
Fix request $i$, the fleet state, the grid, and the same $\theta>0$ in both problems. Write $\boldsymbol q^{(1)}=\widehat{\boldsymbol q}$ and $\boldsymbol q^{(2)}=\boldsymbol q^{\mathbb P}$, with corresponding current fractions $y_i^{(1),*}=\widehat y_i^*$ and $y_i^{(2),*}=y_i^{\mathbb P,*}$. Both plans retain the separate current request with bound one. Use the network in Section~\ref{ec:planning}, including all eligible forecast arcs even when their bounds are zero. For $k\in\{1,2\}$, let $x_{z,\tau|i}^{(k),*}$ be the optimal forecast service amount under forecast $\boldsymbol q^{(k)}$. The contribution of type $z$ in interval $\tau$ to the planning objective \eqref{eq:quadratic-objective} is
\begin{equation}\label{eq:ec-arc-costs}
F_{z,\tau|i}^{(k)}(x):=(r_z+\theta)x -\frac{\theta x^2}{q_{z,\tau|i}^{(k)}},\qquad 0\le x\le q_{z,\tau|i}^{(k)},\quad q_{z,\tau|i}^{(k)}>0.
\end{equation}
If $q_{z,\tau|i}^{(k)}=0$, the service amount is fixed to zero and its objective contribution is zero. Minimizing the negative planning objective gives marginal cost $-\bigl(F_{z,\tau|i}^{(k)}\bigr)'(x)$, where the derivative is with respect to the service amount $x$. Define
\[
\Delta x_{z,\tau|i}:=x_{z,\tau|i}^{(2),*}-x_{z,\tau|i}^{(1),*}, \qquad \Delta q_{z,\tau|i}:=q_{z,\tau|i}^{(2)}-q_{z,\tau|i}^{(1)}.
\]

We first establish a marginal-cost ordering. For $\tau\in\mathcal T_i$ and $z\in\mathcal Z^{\mathrm{end}}_{\tau|i}$, if $|\Delta x_{z,\tau|i}|>|\Delta q_{z,\tau|i}|$, then both forecast bounds are positive and
\begin{equation}\label{eq:ec-marginal-order}
\operatorname{sgn}(\Delta x_{z,\tau|i}) \left[-\bigl(F_{z,\tau|i}^{(2)}\bigr)'(x_{z,\tau|i}^{(2),*}) +\bigl(F_{z,\tau|i}^{(1)}\bigr)'(x_{z,\tau|i}^{(1),*})\right]>0,
\end{equation}
where $\operatorname{sgn}$ denotes the sign function. A sufficiently small change in direction $\operatorname{sgn}(\Delta x_{z,\tau|i})$ at $x_{z,\tau|i}^{(1),*}$, and in the opposite direction at $x_{z,\tau|i}^{(2),*}$, respects the respective arc bounds.

First suppose $\Delta x_{z,\tau|i}>|\Delta q_{z,\tau|i}|$. Feasibility gives
\begin{equation}\label{eq:ec-arc-interiority}
q_{z,\tau|i}^{(1)},q_{z,\tau|i}^{(2)}>0,\qquad x_{z,\tau|i}^{(1),*}<q_{z,\tau|i}^{(1)},\qquad x_{z,\tau|i}^{(2),*}>0.
\end{equation}
Indeed, $q_{z,\tau|i}^{(1)}=0$ would imply $\Delta x_{z,\tau|i}=x_{z,\tau|i}^{(2),*}\le q_{z,\tau|i}^{(2)} =|\Delta q_{z,\tau|i}|$; $q_{z,\tau|i}^{(2)}=0$ would imply $\Delta x_{z,\tau|i}\le0$. Similarly, $x_{z,\tau|i}^{(1),*}=q_{z,\tau|i}^{(1)}$ would imply $\Delta x_{z,\tau|i}\le\Delta q_{z,\tau|i}$. These contradictions establish \eqref{eq:ec-arc-interiority} and the claimed feasible directions. Moreover,
\begin{equation}\label{eq:ec-normalized-flow-order}
\begin{aligned}
q_{z,\tau|i}^{(1)}x_{z,\tau|i}^{(2),*} -q_{z,\tau|i}^{(2)}x_{z,\tau|i}^{(1),*} &=q_{z,\tau|i}^{(1)}\Delta x_{z,\tau|i} -x_{z,\tau|i}^{(1),*}\Delta q_{z,\tau|i}\\
&\ge q_{z,\tau|i}^{(1)} \bigl(\Delta x_{z,\tau|i}-|\Delta q_{z,\tau|i}|\bigr)>0.
\end{aligned}
\end{equation}
Consequently,
\begin{equation}\label{eq:ec-positive-marginal-change}
-\bigl(F_{z,\tau|i}^{(2)}\bigr)'(x_{z,\tau|i}^{(2),*}) +\bigl(F_{z,\tau|i}^{(1)}\bigr)'(x_{z,\tau|i}^{(1),*}) =2\theta\left( \frac{x_{z,\tau|i}^{(2),*}}{q_{z,\tau|i}^{(2)}} -\frac{x_{z,\tau|i}^{(1),*}}{q_{z,\tau|i}^{(1)}}\right)>0.
\end{equation}
Interchanging forecasts 1 and 2 proves the result when $\Delta x_{z,\tau|i}<-|\Delta q_{z,\tau|i}|$.

We now compare flows along cycles, which may contain current-request, forecast, and waiting arcs. To treat these three arc types together, let $\omega$ index the full arc set $\mathcal A_i$, let $\omega_0$ denote the current-request arc, and define the forecast-arc subset $\mathcal A_i^{\mathrm f}:= \{(z,\tau):\tau\in\mathcal T_i,\ z\in\mathcal Z^{\mathrm{end}}_{\tau|i}\}$. Let $\boldsymbol f_i^{(k)}=(f_{\omega|i}^{(k)})_{\omega\in\mathcal A_i}$ be an optimal flow under forecast $k\in\{1,2\}$, with $f_{\omega_0|i}^{(k)}=y_i^{(k),*}$ and $f_{\omega|i}^{(k)}=x_{z,\tau|i}^{(k),*}$ for $\omega=(z,\tau)\in\mathcal A_i^{\mathrm f}$; the remaining components are the optimal waiting flows $w_{s,\tau|i}^{(k),*}$. Write $F_i^{(k)}$ for the planning objective \eqref{eq:quadratic-objective} expressed in these arc flows, and $F_{\omega|i}^{(k)}$ for the contribution of arc $\omega$, so that $F_i^{(k)}(\boldsymbol f_i) =\sum_{\omega\in\mathcal A_i}F_{\omega|i}^{(k)}(f_{\omega|i})$. For the current-request arc, $F_{\omega_0|i}^{(k)}(f)=(r_{\xi_i}+\theta)f-\theta f^2$ for $0\le f\le1$; for a forecast arc $\omega=(z,\tau)\in\mathcal A_i^{\mathrm f}$, $F_{\omega|i}^{(k)}(f)=F_{z,\tau|i}^{(k)}(f)$. Waiting arcs have zero objective contribution and nonnegative flow. We minimize $-F_i^{(k)}$ below and write $\Delta f_{\omega|i}:=f_{\omega|i}^{(2)}-f_{\omega|i}^{(1)}$. The two optimal flows have the same node supplies, and the terminal node absorbs the same total supply. Their difference $\Delta\boldsymbol f_i:=\boldsymbol f_i^{(2)}-\boldsymbol f_i^{(1)}$ therefore has zero net flow at every node. Define the forecast arcs whose nonzero flow changes do not exceed their capacity changes:
\begin{equation}\label{eq:ec-capacity-change-arcs}
\mathcal B_i:= \{(z,\tau)\in\mathcal A_i^{\mathrm f}: 0<|\Delta x_{z,\tau|i}|\le|\Delta q_{z,\tau|i}|\}.
\end{equation}
Every changed forecast arc outside $\mathcal B_i$ satisfies the marginal-cost ordering \eqref{eq:ec-marginal-order} and permits the feasible directions established above. A changed current-request arc also permits small moves toward the other flow at both optima, because its bound is one in both problems, and
\begin{equation}\label{eq:ec-current-marginal-change}
\operatorname{sgn}(\Delta f_{\omega_0|i}) \left[-\bigl(F_{\omega_0|i}^{(2)}\bigr)'(f_{\omega_0|i}^{(2)}) +\bigl(F_{\omega_0|i}^{(1)}\bigr)'(f_{\omega_0|i}^{(1)})\right] =2\theta|\Delta f_{\omega_0|i}|>0.
\end{equation}
For waiting arcs, the same moves respect nonnegativity and the marginal-cost difference is zero.

We next decompose the flow difference into cycles. Since the two plans have the same node supplies, their difference satisfies
\begin{equation}\label{eq:ec-difference-conservation}
\sum_{\omega\text{ entering }v}\Delta f_{\omega|i} =\sum_{\omega\text{ leaving }v}\Delta f_{\omega|i}, \qquad v\in\mathcal N_i.
\end{equation}
Construct an auxiliary network by retaining the direction of each arc with $\Delta f_{\omega|i}>0$, reversing each arc with $\Delta f_{\omega|i}<0$, and omitting zero differences. Assign each retained arc the nonnegative flow $|\Delta f_{\omega|i}|$. Equation~\eqref{eq:ec-difference-conservation} implies that this auxiliary flow has equal inflow and outflow at every node. Thus a walk along positive-flow arcs cannot stop. Because the network is finite, the walk eventually revisits a node. At the first revisit, the segment between the two visits forms a simple directed cycle $Q$. Although the original planning network is acyclic, reversing arcs can create cycles in this auxiliary network. These cycles describe differences between plans, not backward vehicle travel.

For each such cycle, define its signed arc vector $\boldsymbol\chi_i^Q:=(\chi_{\omega|i}^Q)_{\omega\in\mathcal A_i}$, where $\chi_{\omega|i}^Q$ equals $1$ if $Q$ traverses arc $\omega$ in its original direction, $-1$ if it traverses the arc in reverse, and $0$ otherwise. Let $\Delta\boldsymbol f_i^{\mathrm{rem}}$ denote the remaining signed flow difference, initially equal to $\Delta\boldsymbol f_i$. Extract the largest common flow from $Q$:
\begin{equation}\label{eq:ec-cycle-extraction}
\beta_{Q|i}:=\min_{\omega\in Q} |\Delta f_{\omega|i}^{\mathrm{rem}}|>0, \qquad \Delta\boldsymbol f_i^{\mathrm{rem}} \leftarrow\Delta\boldsymbol f_i^{\mathrm{rem}} -\beta_{Q|i}\boldsymbol\chi_i^Q.
\end{equation}
Here $\beta_{Q|i}$ is the amount of flow extracted along cycle $Q$, not an additional decision variable. This subtraction preserves node balance, does not reverse any remaining flow, and makes at least one arc's remaining flow zero. Repeating the procedure therefore exhausts the flow in at most $|\mathcal A_i|$ steps. Let $\mathcal Q_i$ be the resulting finite set of cycles. On each arc, every extracted cycle flow has the sign of $\Delta f_{\omega|i}$, so no cancellation occurs. Consequently,
\begin{equation}\label{eq:ec-cycle-decomposition}
\Delta\boldsymbol f_i =\sum_{Q\in\mathcal Q_i}\beta_{Q|i}\boldsymbol\chi_i^Q,\qquad \beta_{Q|i}>0,\qquad |\Delta f_{\omega|i}|=\sum_{Q\in\mathcal Q_i} \beta_{Q|i}|\chi_{\omega|i}^Q|.
\end{equation}

Suppose a cycle $Q$ contains $\omega_0$ but no arc in $\mathcal B_i$. The feasible directions established above imply that, for a sufficiently small flow amount $\beta_{\mathrm{step}}>0$,
\begin{equation}\label{eq:ec-cycle-feasible-directions}
\boldsymbol f_i^{(1)}+\beta_{\mathrm{step}}\boldsymbol\chi_i^Q \quad\text{and}\quad \boldsymbol f_i^{(2)}-\beta_{\mathrm{step}}\boldsymbol\chi_i^Q
\end{equation}
are feasible under their respective forecasts. Cycle moves preserve all node balances. In particular, no changed zero-bound forecast arc lies on $Q$, since such an arc belongs to $\mathcal B_i$. First-order optimality for the two minimization problems gives
\begin{equation}\label{eq:ec-cycle-optimality}
-\sum_{\omega\in Q}\bigl(F_{\omega|i}^{(1)}\bigr)'(f_{\omega|i}^{(1)}) \chi_{\omega|i}^Q\ge0,\qquad -\sum_{\omega\in Q}\bigl(F_{\omega|i}^{(2)}\bigr)'(f_{\omega|i}^{(2)}) \chi_{\omega|i}^Q\le0.
\end{equation}
Subtracting yields
\begin{equation}\label{eq:ec-cycle-contradiction}
\sum_{\omega\in Q} \left[-\bigl(F_{\omega|i}^{(2)}\bigr)'(f_{\omega|i}^{(2)}) +\bigl(F_{\omega|i}^{(1)}\bigr)'(f_{\omega|i}^{(1)})\right]\chi_{\omega|i}^Q \le0.
\end{equation}
Every forecast term is positive by \eqref{eq:ec-marginal-order}, every waiting term is zero, and the current-request term is strictly positive by \eqref{eq:ec-current-marginal-change}. This contradicts \eqref{eq:ec-cycle-contradiction}. Thus every cycle containing $\omega_0$ meets $\mathcal B_i$.

If $\Delta f_{\omega_0|i}=0$, the claimed bound is immediate. Otherwise, \eqref{eq:ec-cycle-decomposition} and the preceding conclusion give
\begin{equation}\label{eq:ec-fraction-sensitivity}
\begin{aligned}
|y_i^{(2),*}-y_i^{(1),*}| &=|\Delta f_{\omega_0|i}| =\sum_{\substack{Q\in\mathcal Q_i\\\omega_0\in Q}}\beta_{Q|i}\\
&\le\sum_{\substack{Q\in\mathcal Q_i\\\omega_0\in Q}} \beta_{Q|i}\sum_{\omega\in\mathcal B_i}|\chi_{\omega|i}^Q|\\
&\le\sum_{\omega\in\mathcal B_i}\sum_{Q\in\mathcal Q_i} \beta_{Q|i}|\chi_{\omega|i}^Q| =\sum_{\omega\in\mathcal B_i}|\Delta f_{\omega|i}|\\
&\le\sum_{\tau\in\mathcal T_i} \sum_{z\in\mathcal Z^{\mathrm{end}}_{\tau|i}}|\Delta q_{z,\tau|i}| =\varepsilon_i(\boldsymbol q^{(1)},\boldsymbol q^{(2)}).
\end{aligned}
\end{equation}
Finally, $y_i^{(1),*},y_i^{(2),*}\in[0,1]$ also imply $|y_i^{(1),*}-y_i^{(2),*}|\le1$, proving \eqref{eq:calibration-forecast-support}.
\end{proof}

\section{Proof of Theorem~\ref{thm:ts-oracle}}\label{ec:revenue}

\begin{proof}[Proof of Theorem~\ref{thm:ts-oracle}]
We bound the revenue difference by identifying the first request at which the policies make different admission decisions. Give both policies the same historical data and request stream. Write $\widehat a_i$ and $a_i^{\mathbb P}$ for the decisions of $\widehat\pi_{\theta,\gamma}$ and $\pi_{\theta,\gamma}^{\mathbb P}$, respectively, and define their realized revenues by
\begin{equation}\label{eq:ec-realized-revenues}
\widehat W:=\sum_{i=1}^{I}r_{\xi_i}\widehat a_i, \qquad W^{\mathbb P}:=\sum_{i=1}^{I}r_{\xi_i}a_i^{\mathbb P}.
\end{equation}
The event that request $i$ is their first different decision is
\begin{equation}\label{eq:ec-first-difference}
\mathcal K_i^{\mathrm{first}}:=\{i\le I,\ \widehat a_k=a_k^{\mathbb P}\text{ for all }k<i, \ \widehat a_i\ne a_i^{\mathbb P}\}.
\end{equation}
These events are disjoint, so $\sum_{i=1}^{I}\mathbbm{1}\{\mathcal K_i^{\mathrm{first}}\}\le1$. On $\mathcal K_i^{\mathrm{first}}$, identical earlier decisions and the common vehicle assignment rule give the same fleet records before request $i$. The request must be feasible under both policies; otherwise both would reject it.

Let $\widehat y_i^*$ be the planned fraction obtained using the estimated forecast at the reference fleet state. On $\mathcal K_i^{\mathrm{first}}$, the two fractions lie on different sides of the admission threshold. Proposition~\ref{prop:ts-admission} and Lemma~\ref{lem:calibration-forecast-support} therefore give
\begin{equation}\label{eq:ec-first-difference-margin}
|y_i^{\mathbb P,*}-\gamma| \le|\widehat y_i^*-y_i^{\mathbb P,*}| \le\min\left\{1,\varepsilon_i(\widehat{\boldsymbol q},\boldsymbol q^{\mathbb P})\right\}.
\end{equation}
Earlier rewards cancel on $\mathcal K_i^{\mathrm{first}}$, and $0\le r_{\xi_k}\le R$ implies
\begin{equation}\label{eq:ec-remaining-revenue}
|\widehat W-W^{\mathbb P}| =\left|\sum_{k=i}^{I}r_{\xi_k}(\widehat a_k-a_k^{\mathbb P})\right| \le\sum_{k=i}^{I}r_{\xi_k}|\widehat a_k-a_k^{\mathbb P}| \le R(I-i+1).
\end{equation}
If no decision differs, $\widehat W=W^{\mathbb P}$. Combining \eqref{eq:ec-first-difference-margin} and \eqref{eq:ec-remaining-revenue} thus yields, on every realized stream,
\begin{equation}\label{eq:ec-pathwise-revenue}
\begin{aligned}
|\widehat W-W^{\mathbb P}| &\le R\sum_{i=1}^{I}(I-i+1)\mathbbm{1}\{\mathcal K_i^{\mathrm{first}}\}\\
&\le R\sum_{i=1}^{I}(I-i+1)\mathbbm{1}\left\{ |y_i^{\mathbb P,*}-\gamma|\le\min\left\{1,\varepsilon_i(\widehat{\boldsymbol q},\boldsymbol q^{\mathbb P})\right\} \right\}.
\end{aligned}
\end{equation}
Taking expectations and using
\begin{equation}\label{eq:ec-expected-revenue}
\begin{aligned}
\left|V_{\widehat\pi_{\theta,\gamma}}(\mathbb P) -V_{\pi_{\theta,\gamma}^{\mathbb P}}(\mathbb P)\right| &=\left|\mathbb E_{\mathbb P}[\widehat W-W^{\mathbb P}]\right|\\
&\le\mathbb E_{\mathbb P}[|\widehat W-W^{\mathbb P}|]
\end{aligned}
\end{equation}
proves \eqref{eq:ts-oracle-comparison}. The right side of \eqref{eq:ec-pathwise-revenue} is at most $R\sum_{i=1}^{I}(I-i+1)=RI(I+1)/2$, which is integrable under $\mathbb E_{\mathbb P}[I^2]<\infty$.

\end{proof}

\section{Proof of Theorem~\ref{thm:revenue-rate} and tightness}
\label{ec:revenue-convergence}

We prove the historical-data revenue rate in Theorem~\ref{thm:revenue-rate} and establish its tightness.

\subsection{Proof of Theorem~\ref{thm:revenue-rate}}
\label{ec:revenue-rate}

We bound all forecast errors using historical count deviations and the daily request count. Conditional on these quantities, the margin assumption controls the weighted number of decisions that may change. This avoids assuming that online estimation errors and reference margins are independent.

\begin{proof}
The initial fleet and policy parameters are fixed. Under Assumption~\ref{ass:ts-poisson}, independent Poisson increments give the reference policy's true conditional forecasts as
\begin{equation}\label{eq:ec-true-count-forecast}
q^{\mathbb P}_{z,\tau|i} =\int_{t_i+\tau\Delta}^{\min\{H,t_i+(\tau+1)\Delta\}} \alpha_{z,\kappa(u)}\,du.
\end{equation}
These counts exclude the observed current request. The reference policy's trajectory and margins therefore depend only on the operating-day stream, which is independent of the historical counts. This independence is used when conditioning below.

Define the historical contribution to the forecast-error bound by
\begin{equation}\label{eq:ec-historical-forecast-error}
\varepsilon_D^{\mathrm{hist}} := L\Delta\sum_{z\in\mathcal Z}\sum_{j\in\mathcal J} \frac{|N_{z,j}^{\mathrm{hist}}-D h_j\alpha_{z,j}|}{D h_j}.
\end{equation}
Since $N_{z,j}^{\mathrm{hist}}\sim \operatorname{Poisson}(D h_j\alpha_{z,j})$, Cauchy--Schwarz gives
\begin{equation}\label{eq:ec-mean-historical-error}
\begin{aligned}
\mathbb E_{\mathbb P}|N_{z,j}^{\mathrm{hist}}-D h_j\alpha_{z,j}| &\le\sqrt{D h_j\alpha_{z,j}},\\
\mathbb E_{\mathbb P}[\varepsilon_D^{\mathrm{hist}}] &\le\frac{L\Delta}{\sqrt D} \sum_{z\in\mathcal Z}\sum_{j\in\mathcal J} \sqrt{\frac{\alpha_{z,j}}{h_j}}.
\end{aligned}
\end{equation}

The estimator in \eqref{eq:demand-estimate} has denominator at least $D h_j$. Hence, on every realized stream,
\begin{equation}\label{eq:ec-history-online-error}
\begin{aligned}
|\widehat\alpha_{z,j|i}-\alpha_{z,j}| &\le \frac{|N_{z,j}^{\mathrm{hist}}-D h_j\alpha_{z,j}|}{D h_j}\\
&\quad+ \frac{\sum_{k<i}\mathbbm{1}\{\xi_k=z,\ \kappa(t_k)=j\} +\alpha_{z,j}\int_0^{t_i}\mathbbm{1}\{\kappa(u)=j\}\,du} {D h_j}.
\end{aligned}
\end{equation}
Summing the online counts over $z,j$ gives $i-1\le I$, and the corresponding intensity integrals sum to at most $\Lambda$. Every planning window has length at most $L\Delta$, so
\begin{equation}\label{eq:ec-daywise-forecast-error}
\begin{aligned}
\varepsilon_i(\widehat{\boldsymbol q},\boldsymbol q^{\mathbb P}) &\le L\Delta \sum_{z\in\mathcal Z}\sum_{j\in\mathcal J} |\widehat\alpha_{z,j|i}-\alpha_{z,j}|\\
&\le\varepsilon_D^{\mathrm{hist}} +\frac{L\Delta}{D h_{\min}}(I+\Lambda).
\end{aligned}
\end{equation}

Condition on the historical counts and on $I=n\ge1$, and set
\[
\varepsilon_D(n):= \min\left\{1,\varepsilon_D^{\mathrm{hist}} +\frac{L\Delta}{D h_{\min}}(n+\Lambda)\right\}.
\]
This threshold is fixed under the conditioning. Historical independence leaves the conditional distribution of the reference margins unchanged. Assumption~\ref{ass:revenue-margin} and \eqref{eq:ec-daywise-forecast-error} therefore imply
\begin{equation}\label{eq:ec-conditional-margin-bound}
\begin{aligned}
&\mathbb E_{\mathbb P}\!\left[ \sum_{i=1}^{n}(n-i+1) \mathbbm{1}\!\left\{|y_i^{\mathbb P,*}-\gamma|\le \min\{1,\varepsilon_i(\widehat{\boldsymbol q},\boldsymbol q^{\mathbb P})\} \right\} \,\middle|\,(N_{z,j}^{\mathrm{hist}})_{z,j},I=n \right]\\
&\quad\le \mathbb E_{\mathbb P}\!\left[ \sum_{i=1}^{n}(n-i+1)\mathbbm{1}\{|y_i^{\mathbb P,*}-\gamma|\le\varepsilon_D(n)\} \,\middle|\,(N_{z,j}^{\mathrm{hist}})_{z,j},I=n \right]\\
&\quad\le M_{\theta,\gamma}\frac{n(n+1)}{2} \left[\varepsilon_D^{\mathrm{hist}} +\frac{L\Delta}{D h_{\min}}(n+\Lambda)\right].
\end{aligned}
\end{equation}
For $I=0$ both revenues are zero. Taking expectations and applying Theorem~\ref{thm:ts-oracle} yields
\begin{equation}\label{eq:ec-revenue-moment-bound}
\begin{aligned}
&|V_{\widehat\pi_{\theta,\gamma}}(\mathbb P) -V_{\pi_{\theta,\gamma}^{\mathbb P}}(\mathbb P)|\\
&\quad\le R M_{\theta,\gamma} \left[ \mathbb E_{\mathbb P}\!\left[\frac{I(I+1)}{2}\right] \mathbb E_{\mathbb P}[\varepsilon_D^{\mathrm{hist}}] +\frac{L\Delta}{D h_{\min}} \mathbb E_{\mathbb P}\!\left[\frac{I(I+1)(I+\Lambda)}{2}\right] \right].
\end{aligned}
\end{equation}
Here $I$ is independent of $\varepsilon_D^{\mathrm{hist}}$. For $I\sim\operatorname{Poisson}(\Lambda)$,
\begin{equation}\label{eq:ec-poisson-revenue-moments}
\mathbb E_{\mathbb P}\!\left[\frac{I(I+1)}{2}\right] =\frac{\Lambda^2+2\Lambda}{2}, \qquad \mathbb E_{\mathbb P}\!\left[\frac{I(I+1)(I+\Lambda)}{2}\right] =\Lambda^3+3\Lambda^2+\Lambda.
\end{equation}
Substituting \eqref{eq:ec-mean-historical-error} and \eqref{eq:ec-poisson-revenue-moments} into \eqref{eq:ec-revenue-moment-bound} proves \eqref{eq:revenue-learning-rate}.
\end{proof}

\subsection{Tightness of the historical-data revenue rate}
\label{ec:revenue-rate-tightness}

To show that the exponent in Theorem~\ref{thm:revenue-rate} cannot be improved in general, we construct a fixed instance in which the reference policy always accepts its first feasible request, whereas the learned policy can reject it. Both policies earn at most one unit, so losses cannot be canceled by gains on other request streams.

\begin{lemma}[Matching-order revenue lower bound]
\label{lem:revenue-rate-tightness}
There is a fixed instance satisfying Assumptions~\ref{ass:ts-poisson} and~\ref{ass:revenue-margin}, with $\theta=1$ and $\gamma=1/3$, such that
\begin{equation}\label{eq:ec-revenue-lower-limit}
\liminf_{D\to\infty}\sqrt D\, \left[V_{\pi_{1,1/3}^{\mathbb P}}(\mathbb P) -V_{\widehat\pi_{1,1/3}}(\mathbb P)\right] \ge \frac{e^{-3/2}}{2}\bigl[1-\Phi(\sqrt2)\bigr]>0,
\end{equation}
where $\Phi$ is the standard normal distribution function and the learned policy uses $D$ independent complete historical days.
\end{lemma}

\begin{proof}
Take $C=1$, $\Delta=1$, $H=2$, and $L=2$. The vehicle is initially idle in region $s_1$. There is one request type $z$, with $o_z=s_1$, $d_z=s_2\ne s_1$, $\ell_z=1$, and $r_z=1$. There are no requests originating in $s_2$, so at most one ride can be served. Requests after time one cannot finish by $H$. Use three contexts, with the same rates on each independent day:
\begin{equation}\label{eq:ec-lower-contexts}
\kappa(u)=
\begin{cases}
1,&0\le u<1/2,\\
2,&1/2\le u<3/2,\\
3,&3/2\le u<2,
\end{cases}
\qquad (\alpha_{z,1},\alpha_{z,2},\alpha_{z,3})=(1,2,1), \qquad (h_1,h_2,h_3)=(1/2,1,1/2).
\end{equation}
Thus Assumption~\ref{ass:ts-poisson} holds and $\Lambda=3$.

At any feasible request time $0<t_i\le1$, only forecast interval zero is eligible. With $r_z=\theta=1$, the plan uses its full unit of capacity. The parallel-arc identity therefore gives
\begin{equation}\label{eq:ec-lower-fractions}
\begin{aligned}
q_{z,0|i}^{\mathbb P} &=\int_{t_i}^{t_i+1}\alpha_{z,\kappa(u)}\,du =2-|t_i-1/2|,\\
y_i^{\mathbb P,*} &=\frac{1}{3-|t_i-1/2|}\ge\frac13=\gamma, \qquad \widehat y_i^*=\frac{1}{1+\widehat q_{z,0|i}}.
\end{aligned}
\end{equation}
These formulas also hold when an estimated count is zero, in which case the current fraction is one. The reference policy accepts the first request if it arrives by time one; all later requests are infeasible on its trajectory. Consequently,
\begin{equation}\label{eq:ec-lower-pathwise-sign}
0\le\widehat W\le W^{\mathbb P}\le1, \qquad V_{\pi_{1,1/3}^{\mathbb P}}(\mathbb P)=1-e^{-3/2}.
\end{equation}

We next verify Assumption~\ref{ass:revenue-margin}. At a feasible reference decision,
\[
|y_i^{\mathbb P,*}-\gamma|=\frac{|t_i-1/2|}{3(3-|t_i-1/2|)} \ge\frac{|t_i-1/2|}{9}.
\]
Conditional on $I=n$, the $n$ unordered arrival times have density $\alpha_{z,\kappa(u)}/3\le2/3$. If $|y_1^{\mathbb P,*}-\gamma|\le m$, at least one arrival lies within distance $9m$ of $1/2$, an interval of length at most $18m$. A union bound gives $\mathbb P(|y_1^{\mathbb P,*}-\gamma|\le m\mid I=n)\le12nm$. Only the first request can be feasible on the reference trajectory, so
\begin{equation}\label{eq:ec-lower-margin-check}
\mathbb E_{\mathbb P}\!\left[ \sum_{i=1}^{n}(n-i+1)\mathbbm{1}\{|y_i^{\mathbb P,*}-\gamma|\le m\}\,\middle|\,I=n \right] \le12n^2m \le24m\,\frac{n(n+1)}{2}.
\end{equation}
Hence the assumption holds with $M_{1,1/3}=24$, independently of $D$.

The historical context-two count satisfies $N_{z,2}^{\mathrm{hist}}\sim\operatorname{Poisson}(2D)$. Define
\[
\mathcal{K}_D^{\mathrm{hist}}:= \{N_{z,2}^{\mathrm{hist}}\ge2D+2\sqrt D\}, \qquad \mathcal U_D:=\left[\frac12,\frac12+\frac{1}{4\sqrt D}\right].
\]
Let $\mathcal{K}_D^{\mathrm{one}}$ be the event that exactly one request arrives in $[0,1]$ and its time lies in $\mathcal U_D$. The two events are independent. On $\mathcal{K}_D^{\mathrm{one}}$, no operating-day request has been observed before this decision. Estimator~\eqref{eq:demand-estimate} then gives
\begin{equation}\label{eq:ec-lower-overestimate}
\begin{aligned}
\widehat q_{z,0|1} &=(3/2-t_1)\frac{N_{z,2}^{\mathrm{hist}}}{D+t_1-1/2} +(t_1-1/2)\frac{N_{z,3}^{\mathrm{hist}}}{D/2}\\
&\ge \frac{(1-1/(4\sqrt D))(2D+2\sqrt D)}{D+1/(4\sqrt D)} >2 \qquad\text{on }\mathcal{K}_D^{\mathrm{hist}}\cap \mathcal{K}_D^{\mathrm{one}}.
\end{aligned}
\end{equation}
The last inequality holds for $D\ge1$: after subtracting twice the denominator, the numerator difference is $\tfrac32\sqrt D-\tfrac12-\tfrac{1}{2\sqrt D}>0$. Thus $\widehat y_1^*<1/3$ and the learned policy rejects the request, whereas the reference policy accepts it. No later request can be served, giving $W^{\mathbb P}-\widehat W=1$ on this event.

The integrated intensity on $[0,1]$ is $3/2$, and the intensity throughout $\mathcal U_D$ is two. Therefore
\begin{equation}\label{eq:ec-lower-event-probabilities}
\begin{aligned}
\mathbb P(\mathcal{K}_D^{\mathrm{one}})&= \frac{e^{-3/2}}{2\sqrt D},\\
\mathbb P(\mathcal{K}_D^{\mathrm{hist}})&\longrightarrow 1-\Phi(\sqrt2)>0.
\end{aligned}
\end{equation}
The second line follows from the central limit theorem for the sum of the $D$ independent daily Poisson counts. By \eqref{eq:ec-lower-pathwise-sign} and independence,
\[
V_{\pi_{1,1/3}^{\mathbb P}}(\mathbb P) -V_{\widehat\pi_{1,1/3}}(\mathbb P) \ge\mathbb P(\mathcal{K}_D^{\mathrm{hist}})\, \frac{e^{-3/2}}{2\sqrt D}.
\]
This proves \eqref{eq:ec-revenue-lower-limit}. Together with Theorem~\ref{thm:revenue-rate}, it gives a $\Theta(D^{-1/2})$ revenue difference for this fixed instance.
\end{proof}

\section{Computational settings}\label{ec:computation}

\subsection{Implementation and evaluation}
\label{ec:implementation}

Each Gurobi model uses one thread, with barrier for positive regularization and dual simplex for the linear models. The seed is zero and the initial solver tolerances are $10^{-8}$. The solver updates the sparse model at each feasible arrival. In the additional synthetic calibration runs and in New York, a nonoptimal solve is retried after a reset with increased numerical focus. In New York, failure of the planned-fraction or nonboundary decision check also triggers a retry, with barrier tolerance tightened to $10^{-10}$ and then $10^{-12}$ if needed. All accepted solutions must be optimal with primal violation at most $10^{-5}$. The admission tolerance is $10^{-7}$. At positive regularization we check the fraction recovered from $g_i$ against the primal fraction to within $10^{-4}$. Offline solutions must satisfy integrality and flow conservation to within $10^{-6}$. Each simulated pickup is checked against vehicle availability and completion times.

We run up to eight experiments on separate physical cores, varying policy order across streams. Gurobi and BLAS each use one thread per experiment. Timing includes feasibility checks, demand and fleet updates, and all online solves or sampling. We exclude initial model construction and offline calibration. Sequential methods are timed per observed request and Batch per nonempty dispatch epoch. Each reported daily time averages a full operating day.
\subsection{New York service-area construction}
\label{ec:nyc-service-area}

The 25-region area uses taxi-zone identifiers 43, 48, 50, 75, 100, 137, 140, 141, 142, 143, 161, 162, 163, 164, 170, 186, 229, 230, 233, 236, 237, 238, 239, 262, and 263. The geographic construction starts from zone 161 and uses the centroid of the union of zones 161, 162, 236, and 237 as its fixed reference point when selecting adjacent zones. Historical observations determine all type attributes and initial fleet shares. Test-day data enter rate estimates only after their requests have been observed.

\subsection{Benchmark implementation}
\label{ec:benchmark-models}
For Sampled Lookahead, a vehicle enters each sampled network at its first eligible node. Integer supplies and unit request capacities admit an integral optimum. The two actions use the same network with different supplies, giving 40 solves for 20 samples. Sampling uses fitted Poisson rates even for overdispersed demand. The offline bound uses the full observed day and applies to immediate pickups.
\subsubsection{Supply-aware Mirror Backpressure}\label{ec:supply-aware-mbp}
The supply-aware MBP rule of \citet{kanoria2022mirror} uses current idle supply and a price for vehicle occupation. Let $n_{s,i}:=\sum_{c\in\mathcal C}\mathbbm{1}\{s_{c,i}=s,\ u_{c,i}\le t_i\}$ be the idle-vehicle count in region $s$ at arrival $i$. We use
\begin{equation}\label{eq:benchmark-mbp-congestion}
\bar n_{s,i}:=\frac{n_{s,i}+\sqrt C}{(1-\rho)C+S\sqrt C}, \qquad f(x):=-\frac{\sqrt S}{\sqrt x}, \qquad \rho=0.95.
\end{equation}
The congestion difference rewards trips from locations with more idle vehicles to locations with fewer idle vehicles. With step size $\eta=1/C$, the supply price is updated at each arrival by
\begin{equation}\label{eq:benchmark-mbp-price}
v_i^-:=\bigl[v_{i-1}^+-\eta\rho C(t_i-t_{i-1})\bigr]^+, \qquad v_i^+:=v_i^-+\eta\ell_{\xi_i}\Delta a_i,
\end{equation}
where $[x]^+:=\max\{x,0\}$, $t_0=0$, and $v_0^+=0$. We update the supply price at every arrival, including infeasible requests, and apply the admission test in Section~\ref{sec:benchmarks} at feasible arrivals with tolerance $10^{-7}$.
\subsection{Sensitivity to the regularization parameter}
\label{ec:theta-sensitivity}

In a joint calibration, each synthetic family selects $(\theta,\gamma)$ on its two validation streams over $\theta\in\{1,5,10,20,50\}$ and the same 19 thresholds. Ties favor the larger $\gamma$, then $\theta$ closest to 10 in logarithmic distance, then the smaller $\theta$. Four families select the same pair as in the main comparison. Overdispersed selects $(5,0.20)$ and Hotspot $(5,0.30)$, raising mean paired test revenue by 0.11\% and 1.70\%, respectively, over calibrating $\gamma$ at $\theta=10$. The forecasting histories and ten test streams per family are identical across comparisons.

The nominal New York value $\theta=1$ was chosen in a preliminary small-area comparison of $\{1,5,10\}$ with $C=100$ vehicles on January 21--22. It is held fixed in the full-area threshold calibration and test comparisons described above.

\end{document}